\documentclass[english]{article}
\usepackage[T1]{fontenc}
\usepackage[latin9]{inputenc}
\usepackage{units}
\usepackage{textcomp}
\usepackage{ifsym}
\usepackage{amsthm}
\usepackage{amsmath}
\usepackage{amssymb}
\usepackage{graphicx}
\usepackage{setspace}

\makeatletter

\providecommand{\tabularnewline}{\\}

\usepackage{amsmath}
\usepackage{amsfonts}
\usepackage{stmaryrd}
\usepackage{scalerel}

\makeatother

\usepackage{babel}
\begin{document}
\global\long\def\d{\circledast}
 \global\long\def\R{\mathbb{\mathbf{\mathbb{R}}}}
 \global\long\def\i{\mathbf{i}}
 \global\long\def\j{\boldsymbol{j}}
\global\long\def\f{\mathrm{\mathbf{\boldsymbol{\textrm{\ensuremath{\frac{x-y+z}{3}}}}}}}
\global\long\def\p{\frac{\pi}{3}}
\global\long\def\a{\Rightarrow}
\global\long\def\S{\boldsymbol{\boldsymbol{S}}}
 \global\long\def\T{\boldsymbol{\boldsymbol{T}}}
 \global\long\def\X{\boldsymbol{X}}
\global\long\def\Y{\boldsymbol{Y}}
\global\long\def\al{\boldsymbol{\alpha}}
\global\long\def\b{\boldsymbol{\beta}}
\global\long\def\g{\boldsymbol{\boldsymbol{\gamma}}}
\global\long\def\P{\boldsymbol{\boldsymbol{P}}}
\global\long\def\Q{\boldsymbol{\boldsymbol{Q}}}
\global\long\def\bel{\!\in\!}

{\LARGE{}Preface}{\LARGE \par}

The author of the article below, Shlomo Jacobi (1932-2014), passed
away before accomplishing his last mission \textendash{} publishing
his treatise on three-dimensional hypercomplex numbers. Shlomo has
been fascinated with complex numbers since his first year at the Technion
(Israel) in the late 1950s. Although a mechanical engineer by training
and by occupation, he was an inquisitive autodidact, a Renaissance
man who throughout his life accumulated vast knowledge in history,
archeology, physics and other humanistic and exact science fields. 

Following his retirement, Shlomo worked closely with Dr. Michael Shmoish
of the Technion on his mathematical paper during 2012-2013. After
Shlomo's decease (February 2014), Michael took upon himself to complete
Shlomo's unfinished manuscript. 

Out of great respect to Shlomo's passion for math, we \textendash{}
his family \textendash{} felt that his mathematical theory should
be published in a mathematical scientific forum. We are grateful to
Dr. Michael Shmoish for his devotion and belief in this theory, and
for his professional contribution to the publication of this article.
Also we wish to thank Mr. Miel Sharf for his kind help with text preparation.

\newpage{}

\title{On a novel 3D hypercomplex number system}

\author{Shlomo Jacobi}
\maketitle
\begin{abstract}
This manuscript introduces $J_{3}$-numbers, a seemingly missing three-dimensional
intermediate between complex numbers related to points in the Cartesian
coordinate plane and Hamilton's quaternions in the 4D space. The current
development is based on a rotoreflection operator $\,\j\,$ in $\R^{3}$
that induces a novel $\d$-multiplication of triples which turns out
to be \textit{associative}, \textit{distributive} and \textit{commutative.}
\\
This allows one to regard a point in $\R^{3}$ as the three-component
$J_{3}$-number rather than a triple of real numbers. Being equipped
with the $\d$-product, the commutative algebra $\R_{\circledast}^{3}$
is isomorphic to $\R\oplus\boldsymbol{\mathbb{C}}$. Some geometric\textbf{
}{\normalsize{}and} algebraic properties of the $J_{3}$-numbers are
discussed.
\end{abstract}

\section*{Introduction}

It is well-known that any point $\boldsymbol{P}=(u,\, v)\,$ in the
standard Cartesian coordinate system $\{e_{1},\, e_{2}\}\,$ might
be described by a vector $(u,\, v)\,$ or be represented in the form
of a two-term complex number $z=u+\i v$, a scalar. A point $\mathbf{\mathit{\boldsymbol{P}}}=(u,\, v,\, w)\,$
in the standard 3D Cartesian system $\{e_{1},\, e_{2},\, e_{3}\}\,$
is naturally associated with a vector. The question arises: is it
possible to represent the point $\mathbf{\mathit{\boldsymbol{P}}}$
by a \textquotedblleft{}hypercomplex\textquotedblright{} number with
three terms while keeping the basic properties of complex numbers,
including the commutativity of addition and multiplication?

Sir William R. Hamilton, the famous Irish mathematician, for many
years had been trying to extend the algebra of complex numbers to
3D but failed to define the proper multiplication of triples. In 1853
he had published \textquotedblleft{}Lectures on Quaternions\textquotedblright{}
where four-terms \textquotedblleft{}numbers'' were introduced instead
(they were discovered by Hamilton earlier, in 1843, though many believe
that Olinde Rodrigues actually arrived at quaternions without naming
as early as in 1840). Although the quaternions may describe a point
in 4D space and form an associative division algebra over real numbers,
they are not commutative under multiplication, and hence cannot be
really regarded as numbers.

In spite of all failing attempts by Hamilton, we do have an intuitive
conviction that if a point in a plane can be described by a two-term
(complex) number, then \textquotedbl{}there should exist\textquotedblright{}
a three-term-number to represent a point in the Euclidean vector space
$\R^{3}\,$. Thus what we are looking for is a \textquotedbl{}scalar\textquotedbl{}
that fulfills all associative, distributive and commutative laws that
we normally require from numbers. In this paper we present such \textquotedbl{}scalars\textquotedbl{},
coined $J_{3}$-numbers, with an appropriate multiplication rule.
They would imitate complex numbers by keeping a similar pattern with
three components instead of two and having many similar properties.
To describe the $J_{3}$-numbers and their multiplication, a geometric
operator $\j$ acting in $\R^{3}\,$ should be introduced first.

\section{The $\mathbf{\j}$ Operator and \,$J_{3}$-Numbers }

\subsection{Initial definitions}

\medskip{}

Let $\{\boldsymbol{e_{1},\, e_{2},\, e_{3}}\}\,$ be a standard basis
of the real Euclidean vector space $\R^{3}$. Motivated by the action
of the imaginary unit \,$\i$ (regarded as a rotation operator in
the complex plane) we introduce $\j$, a \textit{geometric}\textbf{
operator }acting in $\R^{3}$. It is completely defined by the following
transformation of the basis: 
\[
\boldsymbol{e_{1}\overset{j}{\rightarrow}e_{2}\overset{j}{\rightarrow}e_{3}\overset{j}{\rightarrow}-e_{1}.}\ \tag{1.1}
\]

Let us explain the nature of this linear operator $\j$. First, being
linear it transforms the origin $\boldsymbol{O}=(0,0,0)$$\,$ of
the coordinate system into itself. Second, given a real number $r\text{\ensuremath{\neq}}0$
that describes a distance of $r$ units along the $\boldsymbol{e}_{1}$
axis, the quantity $\j r$ would represent the same distance along
the $\boldsymbol{e}_{2}$ axis, while $\j\j r$ is the distance $r$
on the $\boldsymbol{e_{3}}$ axis. Finally, $\,\j\j\j r\,$ will be
the same distance on the $\boldsymbol{e}_{1}$ axis but in the negative
direction, i.e., $\j\j\j r=-r.\ $ Formally one can conclude that
\[
\j^{3}\overset{def}{=}\j\j\j=-1\ \tag{1.2}
\]
which reminds the famous formula for the imaginary unit: $\,\i^{2}=-1$. 

Now we are ready to introduce the concept of a $\boldsymbol{J_{3}}$\textbf{-number}.

\textbf{Definition 1. }\textit{We will refer to $\S$ of the form
\[
\S=u+\j v+\j\j w,\quad u,\, v,\, w\in\R\ \tag{1.3}
\]
as a $\boldsymbol{J_{3}}$}\textbf{\textit{-number}}\textit{ with
components $\, u,\, v,\, w$. }

\smallskip{}

The following two special $J_{3}$-numbers: $\mathbf{0}$ and $\mathbf{1}$
are given by formulas

\[
\mathbf{0}=0\,+\j0\,+\j\j0,\ \tag{1.3.0}
\]
\[
\mathbf{1}=1\,+\j0\,+\j\j0.\ \tag{1.3.1}
\]

The operator $\j\,$ itself could be also represented by the $J_{3}$-number:

\[
\boldsymbol{\boldsymbol{J}}=\,0\,+\j1\,+\j\j0.\ \tag{1.3.J}
\]

We will regard two $J_{3}$-numbers as equal if their corresponding
components match and define the operations of addition, subtraction
and multiplication by a real scalar component-wise. The relevant commutative,
associative, and distributive laws all follow from the corresponding
laws for the reals. The multiplication of two $J_{3}$-numbers would
be introduced and investigated later.

Since there is a one-to-one correspondence between $J_{3}$-numbers
of the form (1.3), points $\boldsymbol{\boldsymbol{P}}=(u,\, v,\, w)$$\,$
and vectors $\overrightarrow{\P}$ (oriented segments connecting the
origin $\,\boldsymbol{\boldsymbol{O}}\,$ and $\,\boldsymbol{\boldsymbol{P}}$),
from now on we will often use interchangeably words \textquotedbl{}a
point\textquotedbl{}, \textquotedbl{}a vector\textquotedbl{}, \textquotedbl{}a
$J_{3}$-number\textquotedbl{} when dealing with a triple of real
numbers $(u,\, v,\, w)$ and use equalities like $(u,\, v,\, w)=u+\j v+\j\j w$
to emphasize the identity of points/vectors in $\R^{3}\,$ and $J_{3}$-numbers. 

\medskip{}

\textbf{Definition 2.} \textit{The }\textbf{\textit{modulus}}\textit{
of a $J_{3}$-number of the form (1.3) is defined to be the Euclidean
distance from the origin $\boldsymbol{O}\,$ to $\boldsymbol{P}=(u,v,w)$}:

\[
|\S|=\sqrt{u^{2}+v^{2}+w^{2}}.\ \tag{1.4}
\]

\medskip{}

\subsection{Basic properties of the operator $\j$}

\medskip{}

Let us list some basic properties of the$\,\j\,$ operator and \,$J_{3}$-numbers.
It is easily seen from (1.2) and linearity of the operator $\,\j\,$
that

\[
\j\j\j\S=-u-\j v-\j\j w.\ \tag{1.5}
\]

for an arbitrary $J_{3}$-number $\,\S.$ Moreover: 
\[
\j\S=\j(u+\j v+\j\j w)=-w+\j u+\j\j v.\ \tag{1.6}
\]

Much in the same way one can see that 
\[
\j\j\S=-v-\j w+\j\j u.\ \tag{1.7}
\]

Now the following properties of the operator $\,\j\,$ are easy to
prove:

\medskip{}

\textbf{Lemma 1}. \textit{The straight line 
\[
L=\{(x,\, y,\, z):\; x=-y=z\},\ \tag{1.8}
\]
}

\textit{is an invariant 1D subspace of $\R^{3}\,$ under the action
of the operator $\,\j.\,$ Moreover, for any $J_{3}$-number $\,\mathbf{\mathbf{\boldsymbol{\mathbf{\mathit{\boldsymbol{S}}}}}\,}$
on the line $L$$\,$:
\[
\mathbf{\mathit{\j}\S=}\boldsymbol{-}\mathbf{\S.\ \tag{1.9}}
\]
}\textbf{Proof. }Let $\,\S=(s,\,-s,\, s)=s+\j(-s)+\j\j s$, for some
$s$$\,$ in $\R$. Then 

\[
\mathbf{\mathit{\j\S}=}-s+\j s+\j\j(-s)=-(s+\j(-s)+\j\j s)=-\S,
\]

as is easily seen from (1.6). \textit{\scriptsize{}$\text{\textifsymbol[ifgeo]{32}}$}{\scriptsize \par}

\medskip{}

\textbf{Lemma 2}.\textit{ The plane
\[
M=\{(x,\, y,\, z):\; x-y+z=0\},\ \tag{1.10}
\]
}

\textit{is an invariant 2D subspace of $\R^{3}\,$ under the action
of the operator $\,\j.\,$ Moreover, the operator $\,\j\,$ rotates
any nonzero vector in $M\,$ by}\textit{\large{} $\p$}\textit{ radians. }

\textbf{Proof. }Let $\boldsymbol{\mathbf{\mathit{P}}}=p+\j(p+q)+\j\j q\:$
be any \textit{$J_{3}$-number }that belongs to \textit{$M\,$ }for
some real $p$$\:$ and \textbf{$q$}$\,$ and let us denote $\mathbf{\mathit{\boldsymbol{Q}}}=\mathbf{\mathit{\j}\mathbf{\boldsymbol{\mathbf{\mathit{\boldsymbol{P}}}}}}$.
Then 
\[
\boldsymbol{\mathbf{\mathit{\boldsymbol{Q}}}}\overset{(1.6)}{=}-q+\j p+\j\j(p+q)=(-q,\, p,\, p+q)\!\in\! M
\]
as $\,(-q)-p+(p+q)=0$. To complete the proof it is enough to observe
that the triangle $\boldsymbol{\mathbf{\boldsymbol{\mathbf{\mathit{P}}}}\boldsymbol{O}}\mathbf{\mathit{\boldsymbol{Q}}\,}$
is equilateral since both $|\mathbf{\boldsymbol{\mathbf{\mathit{P}}}}|$$\,$
and $|\mathbf{\mathit{\boldsymbol{Q\,}}}|$$\,,$ as well as the distance
between $\mathbf{\boldsymbol{\mathbf{\mathit{\boldsymbol{P}}}}}$
and $\mathbf{\boldsymbol{\mathit{\boldsymbol{Q}}}}$, are all equal
to $\sqrt{(p+q)^{2}+q^{2}+p^{2}}$ (see Fig.1 below). Hence the angle
$\angle\boldsymbol{\mathbf{\boldsymbol{\mathbf{\mathit{P}}}}\boldsymbol{O}}\mathbf{\mathit{\boldsymbol{Q}}=\p}$
. \textit{\scriptsize{}$\text{\textifsymbol[ifgeo]{32}}$}{\scriptsize \par}

\medskip{}

From now on the above subspaces $L$ and $M$ will be called $\j$-invariant
in $\R^{3}.$ 

\medskip{}

Let us observe that the line \textit{$L\,$ }is perpendicular to the
plane \textit{$M\,$ }and that they meet each other at the origin
$\boldsymbol{O.}$ This enables us to interpret the above results
as follows: the operator \textit{$\,\j\,$ }reflects any point on
\textit{$L\,$ }with respect to \textit{$M\,$} and rotates any point
in\textit{ $M\,$} by\textit{\large{} $\p$} radians around the axis
\textit{$L$}. In the general case of an arbitrary point in \textit{$\R^{3}$}
the operator \textit{$\,\j\,$} acts as a \textbf{rotoreflection}.

\smallskip{}

\textbf{Theorem 1.1.}\textit{ The operator $\,\j\,$ reflects any
point in $\R^{3}\,$} \textit{with respect to $M\,$ and subsequently}
\textit{rotates its image by}\textit{\large{} $\p$}\textit{ radians
around the axis $L\,$.}

\textbf{Proof. }Let \textit{$\boldsymbol{X}=(x,\, y,\, z)\,$ }be
any $J_{3}$-number. 

Let us put $\: l=\frac{x-y+z}{3},\; m=\frac{2x+y-z}{3},\; n=\frac{-x+y+2z}{3}\:$
and consider

\[
\boldsymbol{X}_{L}=(l,-l,\, l),\;\boldsymbol{X}_{M}=(m,\, m+n,\, n).\tag{1.11}
\]
One could easily check that $\boldsymbol{X}_{L}\!\in\! L$, $\boldsymbol{X}_{M}\!\in\! M$,
and $\boldsymbol{X}_{L}+\boldsymbol{X}_{M}=\X.$ By linearity of the
operator \textit{$\,\j\,$ }one can write:

\[
\boldsymbol{Y}=\j\boldsymbol{X}=\j\boldsymbol{X}_{L}+\j\boldsymbol{X}_{M}=\boldsymbol{Y}_{L}+\boldsymbol{Y}_{M},
\]

where points $\boldsymbol{Y}_{L},\,\boldsymbol{Y}_{M}\,$ also belong
to \textit{$L\,$ }and \textit{$M$, }respectively, by $\j$-invariance.
To complete the proof it is enough to invoke two previous lemmas (see
also Fig.1 below) .\textit{\scriptsize{} $\text{\textifsymbol[ifgeo]{32}}$}{\scriptsize \par}

\begin{figure}
\includegraphics{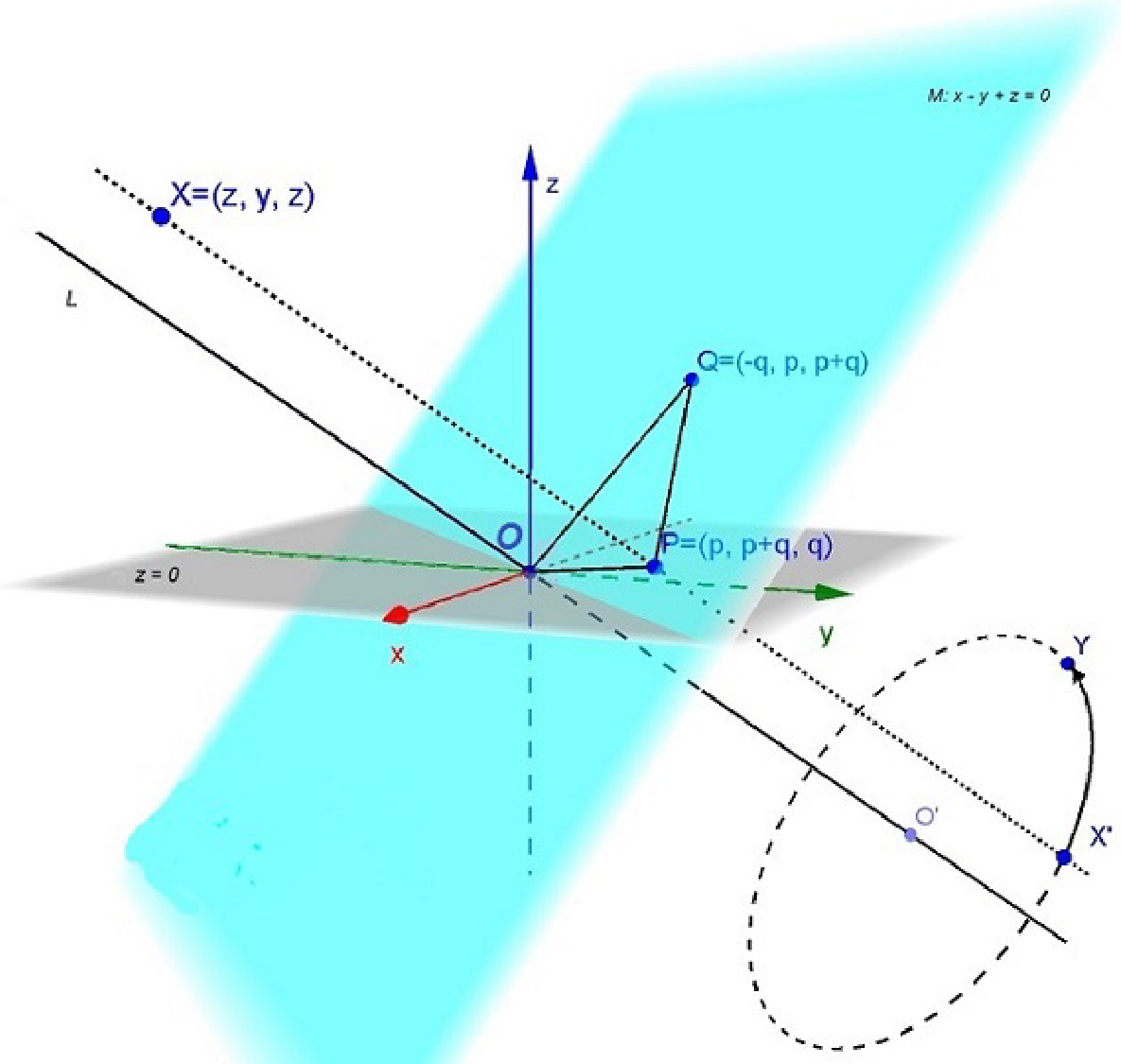}\caption{The operator $\j$ rotates plane $M$ by\textit{\large{} $\p$} radians
since $\boldsymbol{\mathbf{\boldsymbol{\mathbf{\mathit{P}}}}\boldsymbol{O}}\mathbf{\mathit{\boldsymbol{Q}}\,}$
is an equilateral triangle}
\end{figure}

\medskip{}

\subsection{Definition and properties of the \textit{$\d$-}product}

We are ready to define the novel multiplication of $J_{3}$-numbers.
Bearing\textit{ }in mind the multiplication rule for complex numbers
and using (1.2), (1.6), and (1.7) we suggest the following\smallskip{}

\textbf{Definition 3. }\textit{The }\textbf{\textit{$\d$-product}}\textit{
of $J_{3}$-numbers
\[
\boldsymbol{\boldsymbol{T}}=a+\j b+\j\j c,\:\S=u+\j v+\j\j w\ \tag{1.12}
\]
is defined as}
\[
\begin{matrix}\boldsymbol{\boldsymbol{T}}\d\,\S= & (au-bw-cv)+\j(av+bu-cw)+\j\j(aw+bv+cu).\end{matrix}\ \tag{1.13}
\]

\smallskip{}
The right hand side of (1.13) is in the form of a $J_{3}$-number,
thus $\d$-multiplication is a closed operation. Moreover, it is easily
checked to be commutative:

\[
\begin{matrix}\S\d\mathit{\boldsymbol{\boldsymbol{T}}}=\boldsymbol{T}\d\S,\end{matrix}\ \tag{1.14}
\]

and associative:
\[
(\S\d\mathit{\boldsymbol{\boldsymbol{T}}})\d\boldsymbol{U}=\S\d\,(\mathit{\boldsymbol{\boldsymbol{T}}}\d\boldsymbol{U}).\ \tag{1.15}
\]

It is also straightforward to check that all the axioms of the commutative
unital associative algebra over $\R$ hold true for the $\d$-multiplication
and component-wise addition of $J_{3}$-numbers, with zero $\mathbf{0}$
and $\d$-multiplicative unity $\mathbf{1}$ given by formulas (1.3.0)
and (1.3.1), respectively. In particular, for any $\S\in\R^{3}$:

\[
\S+\mathbf{0}=\mathbf{0}+\S=\S,\quad\S+(-\S)=\boldsymbol{0}
\]

\[
\S\d\mathbf{\mathbf{1}}=\mathbf{\mathbf{1}}\d\,\mathit{\S=\S},\ \tag{1.16}
\]

\[
\S\d\mathbf{\mathbf{0}}=\mathbf{\mathbf{0}}\d\,\mathit{\S=\mathbf{\mathbf{0}.}}
\]

\smallskip{}

We will denote this commutative unital associative algebra over $\R$
equipped with the $\d$-product by $\R_{\circledast}^{3}$ or simply
by $\R^{3}$. 

\smallskip{}

Let us also note that any real number $r$ could be uniquely represented
as a $J_{3}$-number with second and third components being equal
to zero:
\[
r\rightarrow\mathbf{r}=(r,\,0,\,0)=r\,+\j0\,+\j\j0,
\]

and that the $\d$-product of such $J_{3}$-numbers by (1.13) reduces
to the usual multiplication of real numbers: $\mathbf{r_{1}\d r_{2}}=r_{1}\, r_{2}.$
In view of (1.3.1) this means that $\R$ is a unital subalgebra of
$\R^{3}$ and the above $r\rightarrow\mathbf{r}$ mapping is a unital
algebra homomorphism. 

\bigskip{}

The $\d$-multiplication table for the standard basis $\{\boldsymbol{e_{1},\, e_{2},\, e_{3}}\}\,$
of the above algebra $\R^{3}$ looks as follows:

\smallskip{}

\begin{doublespace}
\noindent \begin{center}
\begin{tabular}{|c||c|c|c|}
\hline 
$\boldsymbol{\d}$ & $\boldsymbol{e_{1}}$ & $\boldsymbol{e_{2}}$ & $\boldsymbol{e_{3}}$\tabularnewline
\hline 
\hline 
$\boldsymbol{e_{1}}$ & $\boldsymbol{e_{1}}$ & $\boldsymbol{e_{2}}$ & $\boldsymbol{e_{3}}$\tabularnewline
\hline 
$\boldsymbol{e_{2}}$ & $\boldsymbol{e_{2}}$ & $\boldsymbol{e_{3}}$ & \textbf{-}$\boldsymbol{e_{1}}$\tabularnewline
\hline 
$\boldsymbol{e_{3}}$ & $\boldsymbol{e_{3}}$ & \textbf{-}$\boldsymbol{e_{1}}$ & \textbf{-}$\boldsymbol{e_{3}}$\tabularnewline
\hline 
\end{tabular}
\par\end{center}
\end{doublespace}

\[
\mathbf{\mathbf{Table}\:1}
\]

\smallskip{}

The major property of the $\d$-multiplication is given by the following

\textbf{Lemma 3}. \textit{The $\j$-invariant subspaces} \textit{$L\,$
}and \textit{$M$ are mutually $\d$-orthogonal}.

\textbf{Proof. }Let us take arbitrary\textit{ }$J_{3}$-numbers\textit{
}$\boldsymbol{\boldsymbol{L=}}(l,-l,\, l),\;\boldsymbol{\boldsymbol{M=}}(m,\, m+n,\, n)\,$
that belong to \textit{$L\,$ }and \textit{$M,$} respectively. One
can see that \textit{$L\,$ }and \textit{$M$} might be represented
as follows: 
\[
\boldsymbol{\boldsymbol{L=}}l\,(1-\j+\j\j),\,\boldsymbol{M}=(1+\j)\d(m+\j n).
\]

The direct calculations show that
\[
\boldsymbol{\boldsymbol{L}\d}\,\boldsymbol{\boldsymbol{M}}=l\,(1-\j+\j\j)\d(1+\j)\d(m+\j n)=l\,(1+\j^{3})\d(m+\j n)=\boldsymbol{0},
\]

where the last equality is justified by (1.2).\textit{\scriptsize{}
$\text{\textifsymbol[ifgeo]{32}}$}{\scriptsize \par}

\smallskip{}

\subsection{Ideals in $\R^{3}$ }

The above results imply that both \textit{$L\,$ }and \textit{$M$}
are closed under the addition, subtraction and $\d$-multiplication
of $J_{3}$-numbers, i.e., both are subalgebras of $\R^{3}$. Moreover,
they turn out to be ideals of the real algebra $\R^{3}$ due to their
absorbing (or ideal) property:

\[
\forall\S\bel\R^{3},\:\forall\!\boldsymbol{M\bel}M,\:\forall\!\boldsymbol{L\bel}L:\;\S\d\boldsymbol{M\bel}M,\:\S\d\boldsymbol{L\bel}L,
\]

which is obvious in view of \textbf{Lemma 3}. The above inclusions
might be written down in a compact form as follows:
\[
\R^{3}\d M\subset M,\;\R^{3}\d L\subset L.
\]

\smallskip{}

\subsection{Zero divisors and\textit{ $\d$-}invertibility\textit{ }}

\smallskip{}

\textbf{Definition 4.}\textit{ $J_{3}$-number $\boldsymbol{\boldsymbol{S}}$
is called }\textbf{\textit{$\d$-invertible}}\textit{ in }$\R^{3}$\textit{
if there exists another $J_{3}$-number called $\d$-inverse of $\boldsymbol{\boldsymbol{S}}$
and denoted by $\boldsymbol{\boldsymbol{S}}^{\boldsymbol{\mathrm{-1}}}$
such that }

\[
\S^{-1}\d\,\boldsymbol{\boldsymbol{S}}=\boldsymbol{\boldsymbol{S}}\d\,\boldsymbol{\boldsymbol{S}{}^{\mathrm{-1}}}=\boldsymbol{1}.\ \tag{1.17}
\]

\smallskip{}

\textbf{Lemma 4}. \textit{Let $u,\, v,\, w\,$ be real numbers such
that $\Delta=u^{3}-v^{3}+w^{3}+3uvw\neq0$. }

\textit{Then $\:\mathbf{\mathbf{\mathit{\boldsymbol{\boldsymbol{S}}}=}}u+\j v+\j\j w\;$
is $\,\d$-invertible and its $\,\d$-inverse $\boldsymbol{\boldsymbol{S}}^{-1}$
is given by}

\textit{
\[
\boldsymbol{\boldsymbol{S}}^{-1}=\frac{u^{2}+vw}{\Delta}+\j\frac{-w^{2}-uv}{\Delta}+\j\j\frac{v^{2}-uw}{\Delta}.\ \tag{1.18}
\]
}

\textbf{Proof. }Let us rewrite the linear equation $\mathit{\boldsymbol{\boldsymbol{S}}\d\boldsymbol{\! X=}}\boldsymbol{1}$
for an unknown $J_{3}$-number $\boldsymbol{X}\mathbf{\mathbf{=}}\, x+\j y+\j\j z$
in the following form:
\[
(ux-wy-vz)+\j(vx+uy-wz)+\j\j(wx+vy+uz)=1\,+\j0\,+\j\j0.
\]

By matching the corresponding components this might be rewritten as

\[
\left(\begin{matrix}u & -w & -v\\
v & u & -w\\
w & v & u
\end{matrix}\right)\;\left(\begin{matrix}x\\
y\\
z
\end{matrix}\right)=\left(\begin{matrix}1\\
0\\
0
\end{matrix}\right).
\]

To complete the proof it is enough to check that the expression $\Delta=\Delta(u,v,w)$
is a determinant of the Toeplitz matrix 

\[
T=\left(\begin{matrix}u & -w & -v\\
v & u & -w\\
w & v & u
\end{matrix}\right)\ \tag{1.19}
\]
and to invoke the Cramer's rule. \textit{\scriptsize{}$\text{\textifsymbol[ifgeo]{32}}$\smallskip{}
}{\scriptsize \par}

\textbf{Remark. }An intimate connection between $J_{3}$-numbers and
Toeplitz matrices of the form (1.19) will be established below, in
Subsection \textbf{2.2}.

\medskip{}

\textbf{Definition 5. }\textit{A nonzero $J_{3}$-number $\S\neq\mathbf{0}$
is called a }\textbf{\textit{zero divisor}}\textit{ if there exists
$\boldsymbol{\boldsymbol{T}}\neq\mathbf{0}$ in }$\R^{3}$\textit{
such that }

\[
\boldsymbol{\boldsymbol{S}}\boldsymbol{\d}\boldsymbol{\boldsymbol{T}=}\mathbf{0.}\ \tag{1.20}
\]

\medskip{}

We will refer to a pair of non-zero \textit{$J_{3}$}-numbers that
satisfy (1.20) as the \textit{dual} \textit{zero} \textit{divisors}.
The full characterization of zero divisors in\textit{ }$\R^{3}$ and
$\d$-invertible $J_{3}$-numbers is given below.

\textbf{Theorem 1.2}. \textit{Let $\,\mathbf{\mathbf{\mathit{\boldsymbol{\boldsymbol{S}}}=}}u+\j v+\j\j w\neq\boldsymbol{0,}\,$
be a $J_{3}$-number. The following statements are equivalent:}
\begin{enumerate}
\item \textit{The components of $\mathbf{\mathbf{\mathit{\boldsymbol{\boldsymbol{S}}}}}$
satisfy the equation: }
\[
u^{3}-v^{3}+w^{3}+3uvw=0.\ \tag{1.21}
\]

\item $\mathbf{\mathbf{\mathit{\boldsymbol{\boldsymbol{S}}}}}$\textit{
belongs to either one of the $\j$-invariant subspaces} \textit{$L\,$
}and \textit{$M$.}
\item $\mathbf{\mathbf{\mathit{\boldsymbol{\boldsymbol{S}}}}}$ \textit{is
a zero divisor. }
\item $\mathbf{\mathbf{\mathit{\boldsymbol{\boldsymbol{S}}}}}$ \textit{is
not $\d$-invertible}.
\end{enumerate}
\textbf{Proof. }We are going to show that 1$\a$ 2$\a$ 3 $\a$ 4
$\a$1.

1 $\a$2\  The following easy-to-check identity
\[
u^{3}-v^{3}+w^{3}+3uvw=\frac{1}{2}(u-v+w)\,\cdot((u+v)^{2}+(u-w)^{2}+(v+w)^{2})
\]

implies that under (1.21) either $u-v+w=0,$ or $\;(u+v)^{2}+(u-w)^{2}+(v+w)^{2}=0,$
which means that $\,\S$ belongs to either \textit{$M$ }or\textit{
$L$ }as defined by (1.10) or (1.8), respectively\textit{.}

2 $\a$ 3\  It follows from \textbf{Lemma 3} that each point of the
$\j$-invariant subspaces \textit{$L\,$ }and\textit{ $M\,$} ($\,\mathbf{0}$
excluded) is a zero divisor. 

3 $\a$4\  Let us suppose that $\,\S$ is both a zero divisor and
\textit{$\d$}-invertible in $\R^{3}\,$ and take an arbitrary $J_{3}$-number
$\T$ such that $\S\boldsymbol{\d}\,\boldsymbol{\T=}\mathbf{0}$.
Then, in view of (1.15), (1.16), and (1.17), \textit{$\T$} would
also satisfy:

\[
\boldsymbol{\T=(\S{}^{\boldsymbol{\mathrm{-1}}}\d\,\S)}\d\mathit{\T}=\boldsymbol{\S{}^{\boldsymbol{\mathrm{-1}}}\d\:(\S}\d\mathit{\T)=\S{}^{\boldsymbol{\mathrm{-1}}}\d\,}\boldsymbol{0=0.}
\]

Thus, \textit{$\T=0$} which contradicts the assumption that $\,\mathbf{\mathbf{\mathit{\boldsymbol{\boldsymbol{S}}}}}$
is a zero divisor. 

4 $\a$1\  Suppose that (1.21) does not hold. Then \textbf{Lemma
4 }implies that $\mathbf{\mathbf{\mathit{\boldsymbol{\boldsymbol{S}}}}}$
is\textit{ $\d$}-invertible which contradicts the assumption 4. 

This completes the proof of the theorem. \textit{\scriptsize{}$\text{\textifsymbol[ifgeo]{32}}$}{\scriptsize \par}

\medskip{}

\textbf{Lemma }5. \textit{The modulus of the $\d$}-\textit{product
of $\S$ and $\T$ }satisfies the following important inequality:

\[
|\boldsymbol{\boldsymbol{S}}\d\mathit{\boldsymbol{\boldsymbol{T}}}|\leq\sqrt{3}\cdot|\boldsymbol{\boldsymbol{S}}|\,|\mathit{\boldsymbol{\boldsymbol{T}}}|.\ \tag{1.22}
\]

\textbf{Proof.} Let us observe that by direct calculations using (1.4)
and (1.13) one has:

\[
|\boldsymbol{\boldsymbol{S}}\d\mathit{\boldsymbol{\boldsymbol{T}}}|^{2}=|\boldsymbol{\boldsymbol{S}}|^{2}\,|\mathit{\boldsymbol{\boldsymbol{T}}}|^{2}+2\left\{ \S\right\} \left\{ \T\right\} ,
\]
where $\left\{ \mathbf{\mathbf{\mathit{\boldsymbol{\boldsymbol{S}}}}}\right\} =uv-uw+vw,\;\left\{ \boldsymbol{\boldsymbol{T}}\right\} =ab-ac+bc.$
On the other hand, an obvious inequality

\[
(u+v)^{2}+(u-w)^{2}+(v+w)^{2}\geq0
\]

is equivalent to $uv-uw+vw\leq u^{2}+v^{2}+w^{2}$and, therefore,
$\left\{ \mathbf{\mathbf{\mathit{\boldsymbol{\boldsymbol{S}}}}}\right\} \leq\,|\boldsymbol{\boldsymbol{S}}|^{2}$,
$\left\{ \T\right\} \leq\,|\T|^{2}$. Thus, $|\boldsymbol{\boldsymbol{S}}\d\mathit{\boldsymbol{\boldsymbol{T}}}|^{2}\leq3|\boldsymbol{\boldsymbol{S}}|^{2}\,|\mathit{\boldsymbol{\boldsymbol{T}}}|^{2}$
and the rest is plain.\textit{\scriptsize{} $\text{\textifsymbol[ifgeo]{32}}$}{\scriptsize \par}

\smallskip{}

\textbf{Remark}. It follows from \textbf{Lemma 5} that the Hamilton's
\textit{law of moduli} $|\boldsymbol{\boldsymbol{S}}\d\mathit{\boldsymbol{\boldsymbol{T}}}|=|\boldsymbol{\boldsymbol{S}}|\,|\mathit{\boldsymbol{\boldsymbol{T}}}|$
holds true if and only if at least one of $\S$ or $\T$ belongs to
the infinite elliptic cone $xy-xz+yz=0$.

\medskip{}

\textbf{Lemma }6.\textit{ Dual zero divisors belong to different }$\j$\textit{-invariant
subspaces.}

\textbf{Proof. }First, suppose that the dual zero divisors $\boldsymbol{\boldsymbol{L_{1}}}$
and $\boldsymbol{\boldsymbol{L_{2}}}$ belong to the same $\j$-invariant
subspace $L\,$, i.e. $\boldsymbol{\boldsymbol{L_{1}}}=l_{1}(1-\j+\j\j)$
and $\boldsymbol{\boldsymbol{L_{2}}}=l_{2}(1-\j+\j\j)$ for some non-zero
real $l_{1},\, l_{2}.$ It is easy to see that 
\[
\boldsymbol{\boldsymbol{L_{1}}}\d\boldsymbol{\boldsymbol{L_{2}}}=3\, l_{1}\, l_{2}(1-\j+\j\j)=\boldsymbol{0}\ \tag{1.23}
\]

implies that at least one of $l_{1}\,.$ or $l_{2}$ vanishes, which
contradicts the assumption that both $\boldsymbol{\boldsymbol{L_{1}}}$
and $\boldsymbol{\boldsymbol{L_{2}}}$ are zero divisors. 

Now let us take two arbitrary non-zero $J_{3}$-numbers $\boldsymbol{\boldsymbol{M_{1}}}=(m_{1},\, m_{1}+n_{1},\, n_{1})$
and $\boldsymbol{\boldsymbol{M_{2}}}=(m_{2},\, m_{2}+n_{2},\, n_{2})$
from the same $\j$-invariant subspace $M$. 

By \textbf{Lemma 5} the squared modulus of their \textit{$\d$}-product
is equal to:
\[
\boldsymbol{|\boldsymbol{M_{1}}}\d\boldsymbol{\boldsymbol{M_{2}}}|^{2}=\boldsymbol{|\boldsymbol{M_{1}}}|^{2}|\boldsymbol{\boldsymbol{M_{2}}}|^{2}+\boldsymbol{2}\,\left\{ \boldsymbol{\boldsymbol{M_{1}}}\right\} \,\left\{ \boldsymbol{\boldsymbol{M_{2}}}\right\} =\frac{3}{2}\boldsymbol{|\boldsymbol{M_{1}}}|^{2}|\boldsymbol{\boldsymbol{M_{2}}}|^{2},
\]
where the last equality holds true since
\[
\left\{ \boldsymbol{\boldsymbol{M}_{i}}\right\} =m_{i}(m_{i}+n_{i})-m_{i}n_{i}+(m_{i}+n_{i})n_{i}=m_{i}^{2}+m_{i}n_{i}+n_{i}^{2}=\frac{1}{2}\boldsymbol{|\boldsymbol{M}_{i}}|^{2}.
\]
 As $\boldsymbol{\boldsymbol{M_{1}}}$ and $\boldsymbol{\boldsymbol{M_{2}}}$
are assumed to be non-zero their moduli $|\boldsymbol{M_{1}}|$ and
$|\boldsymbol{M_{2}}|$ are both positive. Thus, one has $|\boldsymbol{\boldsymbol{M_{1}}}\d\boldsymbol{\boldsymbol{M_{2}}}|>0$
and consequently $\boldsymbol{\boldsymbol{M}_{1}}\d\boldsymbol{\boldsymbol{M_{2}}}\neq\boldsymbol{0}$,
i.e., $\boldsymbol{\boldsymbol{M_{1}}}$ and $\boldsymbol{\boldsymbol{M_{2}}}$
are not dual zero divisors.\textit{\scriptsize{} $\text{\textifsymbol[ifgeo]{32}}$}{\scriptsize \par}

\newpage{}

\subsection{Linear equations in\textit{ }$\R^{3}$}

\medskip{}

The existence of zero divisors implies that $\R^{3}$ is not a field,
i.e., the division is not always possible. Still, when $\S$ is $\d$-invertible
and $\T$ is an arbitrary $J_{3}$-number then there exists $\X=\S{}^{\boldsymbol{\mathrm{-1}}}\d\,\T$
such that: 

\[
\boldsymbol{\S\d X=\T.}\ \tag{1.24}
\]

We have already characterized in \textbf{Theorem 1.2} all the $\d$-invertible
$J_{3}$-numbers. Now we are going to investigate under what conditions
the linear equation (1.24) is solvable in $\R^{3}.$

\medskip{}

\textbf{Lemma 7}. \textit{Let both $\T=t-\j t+\j\j t\,$ and $\:\S=s-\j s+\j\j s\neq\boldsymbol{0\,}$
belong to the $\j$-invariant subspace} \textit{$L.$ Then the equation
(1.24) has a unique solution in $L$: }

\[
\X_{L}=\frac{t}{3s}\,(1-\j+\j\j)
\]

\textit{and infinitely many solutions in }$\R^{3}:$

\[
\X=\X_{L}+\X_{M},
\]

where $\X_{M}$ is an arbitrary $J_{3}$-number from \textit{the $\j$-invariant
subspace} \textit{$M.$}

\textbf{Proof. }The uniqueness and the formula for \textit{$\X_{L}$}
follows directly from (1.23). The rest is easy due to \textbf{Lemmas
3} and \textbf{6}.\textit{\scriptsize{} $\text{\textifsymbol[ifgeo]{32}}$}{\scriptsize \par}

\medskip{}

\textbf{Lemma 8}. \textit{Let both $\S=u+\j(u+v)+\j\j v\neq\boldsymbol{0}\,$
and $\T=a+\j(a+b)+\j\j b$ belong to the $\j$-invariant subspace}
\textit{$M.$ Then the equation (1.24) has a unique solution in the
plane $M$: }

\[
\X_{M}=c+\j(c+d)+\j\j d,
\]

\[
c=\frac{u\,(2a+b)+v\,(a+2b)}{3(u^{2}+uv+v^{2})},\, d=\frac{u\,(a-b)+v\,(2a+b)}{3(u^{2}+uv+v^{2})}
\]

\textit{and infinitely many solutions in }$\R^{3}:$

\[
\X=\X_{L}+\X_{M},
\]

\textit{where $\X_{L}$ is an arbitrary $J_{3}$-number from} \textit{the
$\j$-invariant subspace} \textit{$L.$}

\smallskip{}

\textbf{Proof. }Let us $\d$-multiply $\S=u+\j(u+v)+\j\j v$ by $\X_{M}=c+\j(c+d)+\j\j d$
and compare the product to \textit{$\T=a+\j(a+b)+\j\j b.$ }By matching
the corresponding coefficients we get a real linear system of 3 equations
with 2 real unknown variables where one equation is a linear combination
of two others. After reducing this system to 2-by-2 case its determinant
happens to be proportional to $u^{2}+uv+v^{2}\neq0.$ Then the formula
for the solution in the plane\textit{ $M$} is a consequence of Cramer's
rule. Now it is enough to invoke \textbf{Lemmas 3} and \textbf{6}
in order to complete the proof.\textit{\scriptsize{} $\text{\textifsymbol[ifgeo]{32}}$}{\scriptsize \par}

\medskip{}

We are going to summarize the above discussion in the following

\textbf{Theorem 1.3}. \textit{Let $\S\neq\boldsymbol{0}$ and $\T$
be }$J_{3}$\textit{-numbers. Then the linear equation }

\[
\boldsymbol{\S\d X=\T.}\ \tag{1.24}
\]

\textit{1) has a unique solution in} $\R^{3}$\textit{ if $\S\,$
is $\d$-invertible;}

\textit{2) has no solution in} $\R^{3}$\textit{ if }

\textit{\qquad{}(a) $\T$ is $\d$-invertible while $\S\,$ is a
zero divisor, or }

\textit{\qquad{}(b) $\S\,$ and $\T$ are dual zero divisors;}

\textit{3) has infinitely many solutions in} $\R^{3}$\textit{ if }

\textit{\qquad{}(a) both $\S\,$ and $\T$ are (non-dual) zero divisors,
such that }$\boldsymbol{\boldsymbol{S}}\d\mathit{\boldsymbol{\boldsymbol{T}}}\neq\boldsymbol{0}$,\textit{
or }

\textit{\qquad{}(b) $\S\,$ is a zero divisor while $\T=\boldsymbol{0}$.}

\textbf{Proof. }1) Given \textit{$\S\,$ }is\textit{ $\d$}-invertible,
let us take a $J_{3}$-number $\X=\S{}^{\boldsymbol{\mathrm{-1}}}\d\,\T.$ 

It follows from (1.16), (1.17), and the associativity of the \textit{$\d$}-multiplication\textit{
}that $\X$ is a solution of (1.24). If $\Y\in\R^{3}$ is also a solution,
i.e., $\S\d\Y=\T$ then

\[
\Y=(\S{}^{\boldsymbol{\mathrm{-1}}}\d\S)\d\Y=\S{}^{\boldsymbol{\mathrm{-1}}}\d(\S\d\Y)=\S{}^{\boldsymbol{\mathrm{-1}}}\d\T=\X.
\]

\textit{2a)} Assume there exists a solution $\X$ of the linear equation
$\boldsymbol{\S\d X=\T}$\textit{ }where \textit{$\T$} is $\d$-invertible
and $\S\,$ is a zero divisor. By setting $\mathit{\boldsymbol{X}}_{1}=\X\d\T^{-1}$one
has:

\[
\S\d\X_{1}=\S\d(\X\d\T^{-1})=(\S\d\X)\d\T^{-1}=\boldsymbol{1},
\]

i.e., \textit{$\S$} turns out to be $\d$-invertible which contradicts
\textbf{Theorem 1.2}. 

\textit{2b)} Now assume that for some $J_{3}$-number $\X$ the equation
(1.24) holds true with $\S\,$ and $\T$ being dual zero divisors. 

Let us observe that by \textbf{Lemma 6} the $J_{3}$-numbers \textit{$\S\,$
}and\textit{ $\T$ }belong to different $\j$-invariant subspaces.
Since \textit{$L\,$ }and \textit{$M$} are ideals in $\R^{3}$ (see
Subsection\textbf{ 1.4}), the left hand side $\boldsymbol{\S\d X}$
belongs to the same $\j$-invariant subspace as \textit{$\S\,$} and
could meet the right hand side \textit{$\T$} of (1.24) only at \textit{$\boldsymbol{0}$}.
This contradicts the assumption that, as a zero divisor,\textit{ }$\boldsymbol{\boldsymbol{T}}\neq\boldsymbol{0}$. 

3) This claim is plain from \textbf{Lemmas 6, 7} and \textbf{8}.\textit{\scriptsize{}
$\text{\textifsymbol[ifgeo]{32}}$}{\scriptsize \par}

\medskip{}

\subsection{Square roots of unity and idempotents in $\R^{3}$}

\medskip{}

Since the $\d$-multiplication is a closed operation in $\R^{3}$
we could define the square of a $J_{3}$-number as the $\d$-product
of $\S=u+\j v+\j\j w\,$ with itself:

\[
\boldsymbol{\boldsymbol{S}^{2}}\overset{def}{=}\boldsymbol{\boldsymbol{S}}\d\mathit{\boldsymbol{\boldsymbol{S}}.}
\]

\medskip{}

\textbf{Theorem 1.4}. \textit{Let $r$ be a real number.} \textit{Then
the quadratic equation in} $\R^{3}$:\textit{ }
\[
\boldsymbol{X^{2}=}r\ \tag{1.25}
\]

\textit{1) has no solution if $r<0$ ,}

\textit{2) has a unique solution $\boldsymbol{X}=\mathbf{0\,}$ if
$\, r=0$, }

\textit{3) has exactly }\textbf{\textit{$4$}}\textit{ distinct solutions
if $\, r>0$.}

\textbf{Proof. }It is easy to check by direct calculations based on
(1.12), (1.13) that 
\[
\boldsymbol{\boldsymbol{S}^{2}}=(u^{2}-2vw)+\j(2uv-w^{2})+\j\j(v^{2}+2uw).\ \tag{1.26}
\]

\smallskip{}

Let us assume that $\boldsymbol{\boldsymbol{S}}\,$ is a solution
of (1.25), i.e, $\boldsymbol{\boldsymbol{S}^{2}}=r+\j0+\j\j0$: 
\[
u^{2}-2vw=r,\;2uv-w^{2}=0,\; v^{2}+2uw=0.\ \tag{1.27}
\]
It follows that:
\[
r=(u^{2}-2vw)-(2uv-w^{2})+(v^{2}+2uw)=(u-v+w)^{2}\geq0.\ \tag{1.28}
\]

This means that the $\d$-square of a $J_{3}$-number cannot be equal
to a negative real number which proves claim 1). 

If $r=0$ the existence of a trivial solution \textit{$\boldsymbol{\boldsymbol{X}}=\mathbf{0\,}$}
is easy. To prove its uniqueness let us assume that there exists an
additional solution $\S\neq\boldsymbol{0}$ of (1.25):

\[
\boldsymbol{\boldsymbol{S}}\d\mathit{\boldsymbol{\boldsymbol{S}}}=\mathbf{0.}
\]

By \textbf{Lemma 6} this implies $\S=\boldsymbol{0}$, as it should
belong to both $\j$-invariant subspaces\textit{ $L\,$ }and\textit{
$M$}. The contradiction proves 2).

Finally, let us consider the remaining case when $\, r>0.$ It is
easy to see from (1.27) that $\, v\,$ and $\, w\,$ are either both
$0,$ or both non-zero. If $v=w=0$ one has $u^{2}=r,\;$ which gives
us two expected real roots: $\pm(\sqrt{r},\;0,\;0).$ If both $\, v\,$
and $\, w\,$ are non-zero then it follows from (1.27) that $v^{2}w^{2}=-4u^{2}vw\Longleftrightarrow vw=-4u^{2}$
while (1.28) shows that $(u-v+w)^{2}=r,$ and the following system
of equations with a positive parameter $r$ emerges: 

\[
9u^{2}=u^{2}-2vw=r,\; vw=-4u^{2}=-\frac{4r}{9},\; v-w=u\pm\sqrt{r}.\ \tag{1.29}
\]

By solving this system one obtains two additional solutions of (1.25)
:

\[
u=\epsilon\,\frac{\sqrt{r}}{3},\; v=\epsilon\,\frac{2\sqrt{r}}{3},\; w=-\epsilon\,\frac{2\sqrt{r}}{3},\;\epsilon=\pm1,\ \tag{1.30}
\]

which completes the proof\textit{\scriptsize{}. $\text{\textifsymbol[ifgeo]{32}}$}{\scriptsize \par}

\textbf{Corollary 1}. There are exactly 4 square roots of unity in
$\R^{3}$:

\[
\boldsymbol{1},\boldsymbol{-1},\,\frac{1}{3}(1+2\j-2\j\j),\,-\frac{1}{3}(1+2\j-2\j\j).\ \tag{1.31}
\]

\textbf{Proof}. The roots $\boldsymbol{1},\boldsymbol{-1}$ are trivial,
the other two are obtained from (1.30) by setting $\, r=1$.\textit{\scriptsize{}
$\text{\textifsymbol[ifgeo]{32}}$}{\scriptsize \par}

\medskip{}

\textbf{Corollary }2. \textit{There are exactly 4 idempotents of algebra}
$\R^{3}$:

\[
\boldsymbol{0},\,\boldsymbol{1},\,\frac{1}{3}(1-\j+\j\j),\,\frac{1}{3}(2+\j-\j\j)\ \tag{1.32}
\]

\textbf{Proof}. Let us use the substitution $\boldsymbol{\X}=\boldsymbol{\Y}+\dfrac{1}{2}$
to rewrite the idempotence equation

\[
\boldsymbol{\X\d\X=}\boldsymbol{\X}\ \tag{1.33}
\]

in the form of (1.25):
\[
\boldsymbol{\Y}^{2}=\frac{1}{4}.
\]

Then it is easily seen that its real roots $\boldsymbol{\Y}=\pm\dfrac{1}{2}$
lead to trivial solutions of (1.33): $\boldsymbol{0}$ and $\boldsymbol{1}$,
while the roots given by (1.30) produce 
\[
\boldsymbol{\al=}\dfrac{1}{3}(1,\,-1,\,1),\;\boldsymbol{\b=}\dfrac{1}{3}(2,\,1,\,-1),\ \tag{1.34}
\]
the two remaining solutions\textit{\scriptsize{}. $\text{\textifsymbol[ifgeo]{32}}$}{\scriptsize \par}

\smallskip{}

Let us note that $\al$ and $\b$ sum up to $\boldsymbol{1}$, belong
to the $\j$-invariant subspaces \textit{$L\,$ }and \textit{$M$},
respectively, and their \textit{$\d$}-product is equal to $\mathbf{0}$:

\[
\boldsymbol{\al=}\dfrac{1}{3}(1-\j+\j\j)\in L,\ \tag{1.35}
\]

\[
\boldsymbol{\b}=\dfrac{1}{3}(2+\j-\j\j)\in M.\ \tag{1.36}
\]

\[
\al\d\boldsymbol{\b=0,\:\al+\b=1.\ }\tag{1.37}
\]

In particular, this leads to two different factorizations of the quadratic 

polynomial $P(\X)=\boldsymbol{X^{2}-}\X$ into linear factors:
\[
P(\X)=\X\d(\X-\boldsymbol{1})=(\X-\al)\d(\X-\b).\ \tag{1.38}
\]

\medskip{}

\subsection{New basis in $\R^{3}$ and the $\d$-multiplication table }

\medskip{}

Let us introduce now the following $J_{3}$-number 
\[
\g=\frac{\sqrt{3}}{3}(\j+\j\j)\in M\ \tag{1.39}
\]

which is orthogonal, as a vector, to both $\al$ and $\b$ and is
$\d$-orthogonal to $\al$ only. It is easy to check that

\[
\b\d\g=\g,\:\g\d\g=-\b\ \tag{1.40}
\]

The triple $\{\boldsymbol{\al},\,\boldsymbol{\b},\,\boldsymbol{\g}\}$
is a new orthogonal (but not $\d$-orthogonal!) basis of $\R^{3}$
and the corresponding $\d$-multiplication table has a particularly
neat form:

\begin{doublespace}
\noindent \begin{center}
\begin{tabular}{|c||c|c|c|}
\hline 
\textbf{$\boldsymbol{\d}$} & \textbf{$\al$} & \textbf{$\b$} & \textbf{$\g$}\tabularnewline
\hline 
\hline 
\textbf{$\al$} & \textbf{$\al$} & \textbf{$\boldsymbol{0}$} & \textbf{$\boldsymbol{0}$}\tabularnewline
\hline 
\textbf{$\b$} & \textbf{$\boldsymbol{0}$} & \textbf{$\b$} & \textbf{$\g$}\tabularnewline
\hline 
\textbf{$\g$} & \textbf{$\boldsymbol{0}$} & \textbf{$\g$} & \textbf{-$\b$}\tabularnewline
\hline 
\end{tabular}
\par\end{center}
\end{doublespace}

\[
\mathbf{\mathbf{Table}\:2}
\]

\smallskip{}

A $J_{3}$-number $\mathbf{\mathbf{\mathit{\boldsymbol{\boldsymbol{S}}}=\,}}u+\j v+\j\j w$
in this basis would be represented as follows:

\[
\S=(u-v+w)\,\al+\frac{1}{2}(2u+v-w)\,\b+\frac{\sqrt{3}}{2}(v+w)\,\g,\ \tag{1.41}
\]

while a $J_{3}$-number $\mathbf{\mathbf{\mathit{\boldsymbol{\boldsymbol{T}}}=\,}}a\cdot\al+b\cdot\b+c\cdot\g$
in the standard basis $\{\boldsymbol{e_{1},\, e_{2},\, e_{3}}\}$
has the following form:

\[
\mathit{\boldsymbol{\boldsymbol{T}}}=\frac{a+2b}{3}+\j\frac{-a+b+c\cdot\sqrt{3}}{3}+\j\j\frac{a-b+c\cdot\sqrt{3}}{3}.\ \tag{1.42}
\]

\smallskip{}

It follows from (1.4) and (1.42) that the modulus of $\mathbf{\mathbf{\mathit{\boldsymbol{\boldsymbol{T}}}}}$
is equal to

\[
|\T|=|a\cdot\al+b\cdot\b+c\cdot\g|=\sqrt{\frac{a^{2}+2b^{2}+2c^{2}}{3}}.\ \tag{1.43}
\]
\[
\]

\subsection{Algebraic connection between $\R^{3}$ and $\mathbb{C}$ \smallskip{}
}

Our next objective is to establish an intimate algebraic connection
between $J_{3}$-numbers and classical real and complex numbers. The
block structure of the Table \textbf{2} hints that we should consider
the direct sum of $\R$ and $\boldsymbol{\mathbb{C}}:$

\[
\mathbf{\mathbf{\mathit{\boldsymbol{D}}}}=\R\oplus\boldsymbol{\mathbb{C}}=\{(r,\:\boldsymbol{z})\,|\, r\in\R,\:\boldsymbol{z}\in\boldsymbol{\mathbb{C}}\}.\ \tag{1.44}
\]

Obviously, it would be a real 3D Euclidean space with addition and
multiplication by a real scalar $k$ defined as

\[
(r_{1},\,\boldsymbol{z}_{1})+(r_{2},\,\boldsymbol{z}_{2})=(r_{1}+r_{2},\,\boldsymbol{z}_{1}+\boldsymbol{z}_{2}),\quad k\,(r,\:\boldsymbol{z})=(k\, r,\: k\,\boldsymbol{z}).\ \tag{1.45}
\]

Let us introduce the product of the above pairs by the following rule:

\[
(r_{1},\,\boldsymbol{z}_{1})\otimes(r_{2},\,\boldsymbol{z}_{2})=(r_{1}\, r_{2},\,\boldsymbol{z}_{1}\,\boldsymbol{z}_{2})\ \tag{1.46}
\]

with $\boldsymbol{z}_{1}\,\boldsymbol{z}_{2}=(x_{1}+\i\, y_{1})\,(x_{2}+\i\, y_{2})\,=(x_{1}x_{2}-y_{1}y_{2})+\i\,(x_{1}y_{2}+x_{2}y_{1})$
being the usual complex number multiplication. It is easily seen that
the product (1.46) is bilinear and thus $\boldsymbol{D}$ becomes
a unital commutative associative algebra over $\R$ with $\boldsymbol{1_{\boldsymbol{D}}=}(1,\,1+\i\,0)$
and $\boldsymbol{0_{\boldsymbol{D}}=}(0,\,0+\i\,0)$ being its unity
and zero, respectively. Note that similarly to $\R^{3}$ the algebra
$\boldsymbol{D}$ also has zero divisors: $(r,\,0)\otimes(0,\boldsymbol{\, z})=\boldsymbol{0}_{\boldsymbol{D}}\;\forall r\in\R,\:\boldsymbol{\forall z}\in\boldsymbol{\mathbb{C}}.$

Now we are ready to establish the isomorphism between $\R^{3}$ and
$\boldsymbol{D}$, all regarded as algebras over $\R$.

\textbf{Theorem 1.5}. $\R^{3}$ is isomorphic to the direct sum of
$\R$ and $\boldsymbol{\mathbb{C}}$:

\[
\R^{3}\cong\R\oplus\boldsymbol{\mathbb{C}.}\ \tag{1.47}
\]

\textbf{Proof}. Let us consider the following mapping $\Phi:\R^{3}\longmapsto\R\oplus\boldsymbol{\mathbb{C}}$
:

\[
\S=u+\j v+\j\j w=r\,\al+x\,\b+y\,\g\longmapsto(r,\, x+\i\, y)=\Phi(\S),\ \tag{1.48}
\]

where the real numbers $r,\, x,\, y$ are the components of $\S$
in the basis $\{\boldsymbol{\al},\,\boldsymbol{\b},\,\boldsymbol{\g}\}$:
\[
r=u-v+w,\, x=\frac{1}{2}(2u+v-w),\, y=\frac{\sqrt{3}}{2}(v+w).\ \tag{1.49}
\]

First of all, formulas (1.41)-(1.42) show that this mapping is a bijection
between $\R^{3}$ and $\boldsymbol{\mathbf{\mathbf{\mathit{\boldsymbol{D}}}}=\R\oplus\boldsymbol{\mathbb{C}}}$.
Secondly, by using (1.48), (1.49) one can see that

\[
\Phi(k_{1}\S+k_{2}\T)=k_{1}\Phi(\S)+k_{2}\Phi(\T),\;\forall k_{1},\, k_{2}\in\R,\;\forall\S,\,\T\in\R^{3},\ \tag{1.50}
\]

and that $\Phi$ maps the unity $\boldsymbol{1}$ of $\R^{3}$ into
the unity of $\boldsymbol{D}$: $\Phi(\boldsymbol{1})=(1,1+\i0)=1_{\boldsymbol{D}}.$ 

Let us observe that by the very definition (1.48):

\[
\Phi(\al)=\Phi(1\al+0\b+0\g)=(1,0),\,\Phi(\b)=(0,1),\,\Phi(\g)=(0,\i).\ \tag{1.51}
\]

It is an easy exercise now to show that the $\d$-multiplication Table
\textbf{2} of Subsection\textbf{ 1.8 }ensures

\[
\Phi(\S\d\T)=\Phi(\S)\otimes\Phi(\T),\;\forall\S,\,\T\in\R^{3},\ \tag{1.52}
\]

which completes the proof.\textit{\scriptsize{} $\text{\textifsymbol[ifgeo]{32}}$}{\scriptsize \par}

\medskip{}

The complex plane $\mathbb{C}$, as an algebra over $\R$, is isomorphic
to the real subalgebra $M$ of $\R^{3}$. This could be seen from
(1.48) where $r$ should be set to $0$ which induces an isomorphism
\[
\Phi_{0}:\: M\longmapsto\mathbb{C}.\ \tag{1.48.0}
\]

Given $\boldsymbol{\boldsymbol{M=}}m+\j(m+n)+\j\j n\in M$, the corresponding
complex number is $\Phi_{0}(\boldsymbol{M})=a+\i b$ with $a=\frac{3}{2}m$
and $b=(\frac{1}{2}m+n)\,\sqrt{3}$. 

Given an arbitrary complex number $\boldsymbol{z=x+\i}y$, the corresponding
$J_{3}$-number in $M$ could be computed by (1.42):
\[
\Phi_{0}^{-1}(x+\i y)=0\,\al+x\,\b+y\,\g=\frac{2x}{3}+\j\frac{x+y\cdot\sqrt{3}}{3}+\j\j\frac{-x+y\cdot\sqrt{3}}{3}.\ \tag{1.53}
\]

In particular, the complex unity $1_{\mathbb{C}}=(1,0)$ corresponds
to the $J_{3}$-number $\b$ which is a unity of the real algebra
$M$, while the imaginary unit $\i=(0,1)$ corresponds to $\g\!\in\!\R^{3}$.

\smallskip{}

We have observed earlier (see Subsection\textbf{ 1.3}) that the algebra
$\R^{3}$ contains $\R$ as a \textbf{\textit{unital}} subalgebra.
Does it similarly contains the real algebra of complex numbers? The
answer is given by the following theorem.

\textbf{Theorem 1.6}. $\mathbb{C}$ is not a \textbf{\textit{unital}}
subalgebra of $\R^{3}.$

\textbf{Proof}. Let us assume that there exists a homomorphism $\Psi:\boldsymbol{\mathbb{C}}\longmapsto\R^{3}$
which maps the complex unity $1_{\mathbb{C}}$ into the unity $\boldsymbol{1}$
of $\R^{3}:$

\[
\Psi(1+\i0)=\boldsymbol{1}.\ \tag{1.54}
\]

Let $\S=\Psi(\i)$ be the $\R^{3}$ image of $\i\in\mathbb{C}$. Then
one has:
\[
\Psi(1_{\mathbb{C}})=\Psi(-\i^{2})=-\Psi(\i)\d\Psi(\i)=-\S\d\S.\ \tag{1.55}
\]

On the other hand, it follows from \textbf{Theorem 1.4} that 
\[
\S^{2}\neq\boldsymbol{-1},\;\forall\S\!\in\!\R^{3}.\ \tag{1.56}
\]

Therefore, $\Psi(1_{\mathbb{C}})$ is not equal to $\boldsymbol{1}$
which contradicts our assumption (1.54). \textifsymbol[ifgeo]{32}

\subsection{Quadratic equations in\textit{ }$\R^{3}$}

We are ready now to investigate the general quadratic equation with
$J_{3}$-coefficients
\[
\boldsymbol{\mathbf{\mathit{A}}}\d\X^{2}+\boldsymbol{\mathit{B}}\d\X+\boldsymbol{C}=\boldsymbol{0}.\ \tag{1.57}
\]

In what follows we are going to use the notion of the $J_{3}$-\textit{discriminant}
which is constructed similarly to the classical case:

\textbf{Definition 6.} \textit{The }$J_{3}$-\textbf{\textit{discriminant}}\textit{
of the equation (1.57) is defined as
\[
\boldsymbol{D}=\boldsymbol{B}^{2}-4\cdot\boldsymbol{A}\d\boldsymbol{C}.\ \tag{1.58}
\]
}

The following notion would also prove useful:

\textbf{Definition 7.} \textit{The }\textbf{\textit{altitude}}\textit{
}$|||\S|||$\textit{ of the $J_{3}$-number $\S=u+\j v+\j\j w$ is
a real number defined as}

\[
|||\S|||=u-v+w.\ \tag{1.59}
\]

\textbf{Remark.} Note that the mapping $\S\rightarrow|||\S|||$ is
a real linear functional on $\R^{3}$. Besides, if $\S=a\cdot\al+b\cdot\b+c\cdot\g$
then by formula (1.49): 

\[
|||\S|||=a.\ \tag{1.60}
\]

Let us also note that since the modulus $|\al|=\frac{\sqrt{3}}{3}$,
the altitude $|||\S|||=a$ is equal to $\sqrt{3}$ times the oriented
distance from the point\textit{ $\S\,$ }to the $\j$-invariant plane
$M$. Obviously, a \textit{$J_{3}$}-number $\S$ has a zero altitude
if and only if $\S\in M$. 

\medskip{}

\textbf{Lemma 9}.\textit{ The altitude of the }\textbf{\textit{$\d$}}\textit{-product
is equal to the product of altitudes:}

\[
|||\S\d\T|||=|||\S|||\cdot|||\T|||,\ \tag{1.61}
\]

\textit{In particular, the altitude of the }\textbf{\textit{$\d$}}\textit{-square
is equal to the square of the altitude:}

\[
|||\boldsymbol{\boldsymbol{S}^{2}}|||\,=|||\S|||^{2}.\ \tag{1.62}
\]

\textbf{Proof. }Since by (1.12), (1.13) one has
\[
\boldsymbol{\boldsymbol{T}}\d\,\S=(au-bw-cv)+\j(av+bu-cw)+\j\j(aw+bv+cu)
\]

it is easy to see that 

\[
|||\S\d\T|||=(a-b+c)(u-v+w)=|||\S|||\cdot|||\T|||,
\]

which proves (1.61) and consequently (1.62).\textit{\scriptsize{}
$\text{\textifsymbol[ifgeo]{32}}$\medskip{}
}{\scriptsize \par}

\textbf{Theorem 1.7. }\textit{The monic equation with} $J_{3}$\textit{-coefficients
}$\P$ \textit{and} $\Q$ :

\[
\X^{2}+\mathit{\P}\d\X+\Q=\boldsymbol{0}\ \tag{1.63}
\]

\textit{is solvable in} \textit{$\R^{3}$ if and only if the altitude
of its }$J_{3}$-\textit{discriminant }$\boldsymbol{D}=\P^{2}-4\cdot\Q$
\textit{is non-negative. More precisely,}

\textit{1) it has a unique solution $\boldsymbol{X}=-\dfrac{1}{2}\P$
if }$\boldsymbol{D}=\boldsymbol{0}$; 

\textit{2) it has exactly }\textbf{\textit{$2$}}\textit{ distinct
solutions if }$\boldsymbol{D}$ \textit{is a zero divisor in} \textit{$\R^{3}$}
\textit{and $\,|||\boldsymbol{D}|||\geq0$};

\textit{3) it has no solutions if $\,|||\boldsymbol{D}|||<0$};\textit{ }

\textit{4) it has exactly }\textbf{\textit{$4$}}\textit{ distinct
solutions if }$\boldsymbol{D}$ \textit{is }$\d$\textit{-invertible
and $\;|||\boldsymbol{D}|||>0$}.

\textbf{Proof. }Let us rewrite (1.63) in the form 
\[
\X^{2}+2\,(\dfrac{1}{2}\,\P)\d\X+(\dfrac{1}{2}\,\P)^{2}=(\dfrac{1}{2}\,\P)^{2}-\Q\boldsymbol{,}\ \tag{1.64}
\]

or, equivalently, after multiplying both sides by $4$:
\[
(2\X+\P)^{2}=\boldsymbol{D}.\ \tag{1.65}
\]

If (1.63) is solvable and $\X$ is a root, then in view of \textbf{Lemma
9} one has:

\[
|||\boldsymbol{D}|||=|||(2\X+\P)^{2}|||=|||2\X+\P|||^{2}\geq\boldsymbol{0}\ \tag{1.66}
\]

which proves the necessary condition and, equivalently, the claim
3). 

To prove the sufficiency it is enough to justify the remaining three
claims. 

1) If $\P$ and $\Q$ are such that $\boldsymbol{D}=\boldsymbol{0}$
then in view of\textbf{ Theorem 1.4} one has $2\X+\P=\boldsymbol{0}$
, i.e., \textit{$\boldsymbol{X}=-\dfrac{1}{2}\P$ }is a unique solution
of (1.63).

2) Assume now that $\boldsymbol{D}=d_{1}\al+d_{2}\b+d_{3}\g\neq\boldsymbol{0,}\:|||\boldsymbol{D}|||\geq0$
and let us write down the unknown $J_{3}$-number $2\X+\P$ in the
basis $\{\boldsymbol{\al},\,\boldsymbol{\b},\,\boldsymbol{\g}\}$:

\[
2\X+\boldsymbol{P}=a\cdot\al+b\cdot\b+c\cdot\g,\ \tag{1.67}
\]

where $a,\, b,\, c$ are real unknowns. By plugging $\boldsymbol{D}$
and $2\X+\P$ into the equation (1.65) one gets:

\[
(a\cdot\al+b\cdot\b+c\cdot\g)^{2}=d_{1}\al+d_{2}\b+d_{3}\g,\ \tag{1.68}
\]

and after invoking the $\d$-multiplication Table \textbf{2} the following
equation emerges:

\[
a^{2}\al+(b^{2}-c^{2})\b+2bc\g=d_{1}\al+d_{2}\b+d_{3}\g.\ \tag{1.69}
\]

If $\boldsymbol{D}$ is a zero divisor then by \textbf{Theorem 1.2}
it belongs either to the $\j$-invariant subspace\textbf{ $M$ }or
to\textbf{ $L$.} 

We will consider each case separately:

a) Assume that $\boldsymbol{D\!\in\!}M$, then one has $d_{1}=|||\boldsymbol{D}|||=0$
which implies $a=0$ and the equation (1.69) could be rewritten as:

\[
(b\cdot\b+c\cdot\g)^{2}=d_{2}\b+d_{3}\g.\ \tag{1.70}
\]

Let us denote $\boldsymbol{z}=\Phi_{0}(b\cdot\b+c\cdot\g)=b+\i c\in\!\mathbb{C}$
and apply the isomorphism (1.48.0) to both sides of (1.70). Then we
get the following simple equation in complex numbers:
\[
z^{2}=d_{2}+\i d_{3},\ \tag{1.71}
\]
 where at least one of the components $d_{2},\, d_{3}$ is nonzero
since $\boldsymbol{D}\neq\boldsymbol{0}.$ Of course, this equation
has 2 distinct complex roots $z_{1},\, z_{2}$ and consequently (1.63)
also has exactly 2 distinct solutions:
\[
\boldsymbol{Z}_{i}=\frac{1}{2}(-\boldsymbol{P}+\Phi_{0}^{-1}(z_{i}))\!\in\! M,\: i=1,2.\ \tag{1.72}
\]

b) Assume that $\boldsymbol{D\!\in\!}L$ and \textit{$|||\boldsymbol{D}|||>0$}.
In this case $d_{1}>0$ while $d_{2}=d_{3}=0$ and the equation (1.69)
could be rewritten as:

\[
a^{2}\al+(b^{2}-c^{2})\b+2bc\g=d_{1}\al+0\b+0\g,\ \tag{1.73}
\]

which obviously has 2 distinct solutions: $\pm\sqrt{d_{1}}\al+0\b+0\g.$ 

This proves the claim 2). 

To prove 4) it is enough to observe that if $d_{1}=|||\boldsymbol{D}|||>0$
and $\boldsymbol{D}$ is $\d$\textit{-invertible }then at least one
of the components $d_{2},\, d_{3}$ is nonzero, and the equation (1.69)
breaks down into two independent equations

\[
a^{2}\al=d_{1}\al,\;(b\cdot\b+c\cdot\g)^{2}=d_{2}\b+d_{3}\g\ \tag{1.74}
\]
each one with exactly two distinct solutions. Consequently, there
are exactly four distinct solutions of the equation (1.69), namely:

\[
\X_{1}=-\sqrt{d_{1}}\al+\boldsymbol{Z}_{1},\,\X_{2}=\sqrt{d_{1}}\al+\boldsymbol{Z}_{1},
\]
\[
\X_{3}=-\sqrt{d_{1}}\al+\boldsymbol{Z}_{2},\,\X_{4}=\sqrt{d_{1}}\al+\boldsymbol{Z}_{2},\ \tag{1.75}
\]

where $\boldsymbol{Z}_{1},\,\boldsymbol{Z}_{2}$ are given by (1.72).\textit{\scriptsize{}
$\text{\textifsymbol[ifgeo]{32}}$}{\scriptsize \par}

\smallskip{}

Let us discuss the general quadratic equation (1.57). First, we will
note that when its leading coefficient $\boldsymbol{A}$ is $\,\d$-invertible
one can write down the equivalent monic equation:
\[
\X^{2}+\mathit{\P'}\d\X+\Q'=\boldsymbol{0,}\ \tag{1.76}
\]
where $\P'=\boldsymbol{A}^{-1}\d\boldsymbol{B},\:\Q'=\boldsymbol{A}^{-1}\d\mathbf{C}$
and apply the \textbf{Theorem 1.7}.

\smallskip{}

If $\boldsymbol{A}$ is a zero divisor, however, the investigation
of the root structure of (1.57) becomes more involved. 

We will outline the possible research direction in the particular
case when both $\boldsymbol{A}$ as well as $\boldsymbol{B}$ belong
to the subalgebra $M\bel\R^{3}$. Under these conditions our equation
(1.57) could be written in the following form:

\[
(\boldsymbol{A\d}\X+\boldsymbol{B})\d\X=-\boldsymbol{C},\ \tag{1.77}
\]

where by the ideal property of $M$ one has $(\boldsymbol{A\d}\X+\boldsymbol{B})\!\in\! M,\;\forall\X\!\in\!\R^{3}$. 

The following three cases are to be considered:

1) If\textit{ }$\boldsymbol{C}\!\notin\! M$ then there is no solution
since $(\boldsymbol{A\d}\X+\boldsymbol{B})\d\X\!\in\! M,\;\forall\X\!\in\!\R^{3}$.

2) If $\,\boldsymbol{C}=\boldsymbol{0\:}$ then the above equation
becomes:
\[
(\boldsymbol{A\d}\X+\boldsymbol{B})\d\X=\boldsymbol{0}.\ \tag{1.77.0}
\]

By \textbf{Lemma 8} the linear equation $\boldsymbol{A\d}\X=-\boldsymbol{B}$
has infinitely many solutions that fill the straight line

\[
\X_{0}+L=\left\{ \X_{0}+\boldsymbol{L}\,|\;\boldsymbol{L}=l\,(1-\j+\j\j),\: l\!\in\!\R\right\} 
\]

where $\X_{0}$ is its unique solution that belongs to $M$. 

In addition, since $(\boldsymbol{A\d}\X+\boldsymbol{B})\!\in\! M$
the whole line $L$ consists of solutions of (1.77.0) due to \textbf{Lemma
3}. These two solution lines coincide if $\X_{0}=\boldsymbol{0}$,
which is the case when $\boldsymbol{B}=\boldsymbol{0}$.

3) If $\boldsymbol{C}\!\in\! M$ and $\,\boldsymbol{C}\neq\boldsymbol{0}$,
the coefficients of (1.57) are factorized as follows:

\[
\boldsymbol{A=}(1+\j)\d(a_{1}+\j a_{2}),\:\boldsymbol{B=}(1+\j)\d(b_{1}+\j b_{2}),\,\boldsymbol{C=}(1+\j)\d(c_{1}+\j c_{2}),\ \tag{1.78}
\]

and our equation (1.57) becomes:
\[
(1+\j)\d(\boldsymbol{\mathbf{\mathit{A}}}'\d\X^{2}+\boldsymbol{\mathit{B}}'\d\X+\boldsymbol{C}')=\boldsymbol{0},\ \tag{1.79}
\]

where $\boldsymbol{\mathbf{\mathit{A}}}'=a_{1}+\j a_{2},\,\boldsymbol{\mathit{B}}'=b_{1}+\j b_{2},$
and $C'=c_{1}+\j c_{2}\neq\boldsymbol{0}$. The equality (1.79) tells
us that for any fixed $\X\!\in\!\R^{3}$ the \textit{$J_{3}$-}number
$F(\X)=\boldsymbol{\mathbf{\mathit{A}}}'\d\X^{2}+\boldsymbol{\mathit{B}}'\d\X+\boldsymbol{C}'$
is $\,\d$-orthogonal to $\boldsymbol{I=}1+\j\in\! M$ and thus, by
Lemma \textbf{6}, it should necessarily belong to $L$. Therefore
the above equation (1.79) turns out to be equivalent to 

\[
\boldsymbol{\mathbf{\mathit{A}}}'\d\X^{2}+\boldsymbol{\mathit{B}}'\d\X+\boldsymbol{C}'=\boldsymbol{L},\ \tag{1.80}
\]

where $\boldsymbol{L}=l\,(1-\j+\j\j)$ is an arbitrary $J_{3}$-number
in $L$. 

There are two cases to work out:

a) When $a_{1}\neq a_{2}$ the leading coefficient $\boldsymbol{\mathbf{\mathit{A}}}'=a_{1}+\j a_{2}\!$
is clearly $\,\d$-invertible, and one could reduce (1.80) to the
monic equation

\[
\X^{2}+\mathit{\P}''\d\X+\Q''=\boldsymbol{0},\ \tag{1.81}
\]

with $\P''=\boldsymbol{A}'^{-1}\d\boldsymbol{B}',\:\Q''=\boldsymbol{A}'^{-1}\d(\mathbf{C}'-\boldsymbol{L})$
where $\boldsymbol{L\!\in\!}L$ is a free parameter, and apply the
\textbf{Theorem 1.7}. 

b) The remaining case when $a_{1}=a_{2}$, i.e., $\boldsymbol{\mathbf{\mathit{A}}}'\!\in\! M$
could be treated after writing down the coefficients in the basis
$\{\al,\,\b,\:\g\}$ similarly to the proof of the \textbf{Theorem
1.7}. The details are left to the reader.

\medskip{}

\section{Matrix Representation and Conjugates.}

We will start this subsection by reminding few basic facts from the
theory of the complex numbers.

\subsection{Basic facts on complex numbers}

A complex number $z=a+b\i$ can be represented by the following $2\times2$
matrix:

\[
Z=\left(\begin{array}{cc}
a & -b\\
b & a
\end{array}\right),\ \tag{2.1}
\]

while the conjugate $\bar{z}=a-b\i$ corresponds to the transpose
of the above matrix:

\[
Z^{t}=\left(\begin{array}{cc}
a & b\\
-b & a
\end{array}\right)\ \tag{2.1t}
\]

and is a reflection of the point $z$ across the real axis in the
complex plane.

Let us also remind the following important property of the complex
numbers:

\[
z\bar{z}=|z|^{2}=det\,(Z)\:\in\R,\ \tag{2.2}
\]

where $|z|=\sqrt{a^{2}+b^{2}}$ is an absolute value of $z$ and $det$
stands for the determinant of a matrix.

\subsection{Matrix representation in $\R^{3}$
\[
\]
}

Now we are going to introduce a similar matrix representation of \textit{$J_{3}$-}numbers.

We argue that the Toeplitz matrices of the following special structure

\[
T=\left(\begin{matrix}u & -w & -v\\
v & u & -w\\
w & v & u
\end{matrix}\right),\ \tag{2.3}
\]

which appeared in the proof of \textbf{Lemma 4} (see Subsection \textbf{1.5}),
correspond to the \textit{$J_{3}$-}numbers very much like the matrices
of the form (2.1) are related to complex numbers.

First of all, it is easy to prove that the set of real Toeplitz matrices
of the form (2.3) contains an identity matrix ($u=1,\, v=w=0$) and
is closed under the standard matrix addition and multiplication. Thus
it is a subalgebra of the algebra of real matrices $\R^{3\times3}$.
Secondly, the bijective map
\[
\R^{3}\ni\boldsymbol{\boldsymbol{T}}=u+v\j+w\j\j\leftrightsquigarrow T=\left(\begin{matrix}u & -w & -v\\
v & u & -w\\
w & v & u
\end{matrix}\right)\in\R^{3\times3}\ \tag{2.4}
\]

is a unital algebra isomorphism since the matrix $Q\,$ related to
the $J_{3}$-number $\boldsymbol{\boldsymbol{Q=T}}\d\,\mathit{\boldsymbol{\boldsymbol{S}}}\,$
is equal to the product of matrices $T\,$ and $S\,$ of the form
(2.3) as can be easily checked. In particular, the multiplication
of the Toeplitz matrices of the form (2.3) turns out to be commutative
and any integer power of $\T$ corresponds to the same power of the
matrix $T$: 
\[
\boldsymbol{\boldsymbol{T}}^{n}\leftrightsquigarrow T^{n},\, n=1,\,2,\ldots\ \tag{2.5}
\]

Most of the facts established in Section \textbf{1} for the $J_{3}$-numbers
could be now reformulated in terms of the above Toeplitz matrices.
For example,\textbf{ Lemma 1 }means that the matrix
\[
J=\left(\begin{matrix}0 & 0 & -1\\
1 & 0 & 0\\
0 & 1 & 0
\end{matrix}\right)\ \tag{2.6}
\]

which is related to our basic operator $\j$, has an eigenvalue $-1$
corresponding to the eigenvector $\al=\frac{1}{3}(1,-1,1)$. In addition,
the absorbing property of the $\j$-invariant subspace $L$ (see Subsection
\textbf{1.4}) simply means that $\al$ is an eigenvector of \textbf{any
}Toeplitz matrix $T$ of the form (2.3) with an eigenvalue equal to
$u-v+w$, the altitude of $\T$: $\al\, T=(u-v+w)\,\al.$ 

On the other hand, the matrix theory might help to treat the $J_{3}$-numbers.
It seems natural to define the conjugate of a\textit{ $J_{3}$-}number
$\boldsymbol{\boldsymbol{T}}$ via the transpose of the corresponding
Toeplitz matrix $T$:

\[
T^{t}=\left(\begin{matrix}u & v & w\\
-w & u & v\\
-v & -w & u
\end{matrix}\right).\ \tag{2.7}
\]

\subsection{Conjugation in $\R^{3}$ and its geometric meaning}

\smallskip{}

\textbf{Definition 2.1. }\textit{The $J_{3}$}\textbf{\textit{-conjugate}}\textit{
of $\,\T=u+\j v+\j\j w\:\in\R^{3}$ is defined as:}

\[
\T^{*}=u-\j w-\j\j v.\ \tag{2.8}
\]
\smallskip{}

Let us observe that in view of (1.42) the $J_{3}$-numbers $\boldsymbol{\al},\,\boldsymbol{\b},\,\boldsymbol{\g}$
which constitute the basis $\{\boldsymbol{\al},\,\boldsymbol{\b},\,\boldsymbol{\g}\}$
have simple \textit{$J_{3}$-}conjugates:

\[
\al^{*}=\al,\,\b^{*}=\b,\,\g^{*}=-\g.\ \tag{2.9}
\]

Thus, the \textit{$J_{3}$-}conjugate of $\S=a\cdot\al+b\cdot\b+c\cdot\g$
has the following form:

\[
\S^{*}=a\cdot\al+b\cdot\b-c\cdot\g,\ \tag{2.10}
\]

and obviously the \textit{$J_{3}$-}conjugation is an involution,
i.e., the \textit{$J_{3}$-}conjugate of $\S^{*}$ is equal to $\S$,
which resembles the classical case: $\bar{\bar{z}}=z$ for $z\in\mathbb{C}$
and allows us to talk about \textit{$J_{3}$-}conjugate pairs.

The geometric meaning of the \textit{$J_{3}$-}conjugates is becoming
clear. Namely, since $\{\boldsymbol{\al},\,\boldsymbol{\b},\,\boldsymbol{\g}\}$
is an orthogonal basis in $\R^{3}$ the \textit{$J_{3}$-}conjugation
leads to a reflection through the plane $R$ which spans vectors $\boldsymbol{\al},$$\,\b$
(and $\boldsymbol{e_{1}}=\al+\b=\boldsymbol{1}$). Having $\g=\dfrac{1}{\sqrt{3}}(0,\,1,\,1)$
as a normal vector, the plane $R$ could be described in the original
3D Cartesian coordinate system $Oxyz$ by the following equation:

\[
y+z=0.\ \tag{2.11}
\]

We will refer to this plane $R$ as the \textbf{\textit{conjugate}}
plane.

Note that due to well-known properties of matrix transposes the \textit{$J_{3}$-}conjugation
distributes over the addition and $\d$-multiplication:

\[
(\S\pm\T)^{*}=\S^{*}\pm\T^{*},\,(\S\d\T)^{*}=\S^{*}\d\T^{*}.\ \tag{2.12}
\]

Moreover, it follows from (2.9) that $\S=\S^{*}$ if and only if $\S\!\in\! R$. 

In addition, let us observe that 
\[
\T\d\T^{*}=(u^{2}+v^{2}+w^{2})+(\j-\j\j)(uv-uw+vw)=|\T|^{2}+(\j-\j\j)\{\T\}.\ \tag{2.13}
\]

Due to (2.11) this means that the \textit{$\d$-}product of any \textit{$J_{3}$-}conjugate
pair belongs to the conjugate plane $R$.

\smallskip{}

Alternative way to see that is to consider $\S=a\cdot\al+b\cdot\b+c\cdot\g,\:\S^{*}=a\cdot\al+b\cdot\b-c\cdot\g$
and to invoke the multiplication Table \textbf{2}:

\[
\S\d\S^{*}=a^{2}\al+(b^{2}+c^{2})\b\in R.\ \tag{2.14}
\]

\textbf{Remark} \textbf{1}. By comparing (2.2) and (2.14), one can
regard the plane $R\!\in\!\R^{3}$ given by (2.11) as a counterpart
of the real axis in $\mathbb{C}$.

\textbf{Remark 2.} The following equality holds true:

\[
|\boldsymbol{\mathbf{\mathit{\boldsymbol{\boldsymbol{T}}}}}|^{2}=\boldsymbol{\boldsymbol{\boldsymbol{T}}}\d\,\T^{*},\ \tag{2.15}
\]

if and only if $\boldsymbol{\boldsymbol{T}}$ belongs to the elliptic
cone 
\[
xy-xz+yz=0.\ \tag{2.16}
\]

\textbf{Remark} \textbf{3}. The altitude of the \textit{$\d$-}product
of any \textit{$J_{3}$-}conjugate pair $\S,\:\S^{*}$ is non-negative:
$|||\S\d\S^{*}|||\geq0$$ $ as easily seen from (1.60) and (2.14).

\smallskip{}

\subsection{Determinant of $J_{3}$-numbers }

\smallskip{}

In addition to the notions of the modulus and the altitude of a $J_{3}$-number
let us also introduce the determinant $||\cdot||$.

\smallskip{}

\textbf{Definition 2.2.} The \textbf{determinant} $||\T||$ of a $J_{3}$-number
$\boldsymbol{\boldsymbol{T}}=u+v\j+w\j\j$ is defined to be a determinant
of the corresponding Toeplitz matrix:

\[
||\T||=det\,(T).\ \tag{2.17}
\]

Note that $||\T||\neq0$ if and only if $\T$ is an invertible $J_{3}$-number
(see Theorem \textbf{1.2}). Besides, due to the well-known property
of the matrix determinants one has:
\[
||\T\d\S||=det\,(T\, S)=det(T)\, det(S)=||\T||\cdot||\S||\ \tag{2.18}.
\]

\smallskip{}

Our next objective is to express the above determinant in terms of
the components $u,\, v,\, w$ of a $J_{3}$-number or, alternatively,
in terms of its coefficients $a,\, b,\, c$ in the basis $\{\boldsymbol{\al},\,\boldsymbol{\b},\,\boldsymbol{\g}\}$.

\textbf{Lemma 10}. Let $\boldsymbol{\boldsymbol{T}}=u+v\j+w\j\j=a\cdot\al+b\cdot\b+c\cdot\g$
be a $J_{3}$-number. Then its determinant can be computed as follows:

\[
||\T||=u^{3}-v^{3}+w^{3}+3uvw,\ \tag{2.19}
\]
\[
||\T||=a\cdot(b^{2}+c^{2}).\ \tag{2.20}
\]

\textbf{Proof}. The corresponding Toeplitz matrix $T$ of the form
(2.3) has a determinant which is equal to the right hand side of (2.19)
and so does $||\T||$. In order to prove (2.20) it is enough to express
$a,\, b,\, c\;$ in terms of $u,\, v,\, w\;$ according to (1.41)
and simplify .\textit{\scriptsize{} $\text{\textifsymbol[ifgeo]{32}}$}{\scriptsize \par}

\smallskip{}

Note that it is more instructive to prove (2.20) by figuring out the
following simple structure of an image $D$ of the above Toeplitz
$T$ :

\[
D=\left(\begin{matrix}a & 0 & 0\\
0 & b & -c\\
0 & c & b
\end{matrix}\right)\ \tag{2.21}
\]

under the determinant-preserving linear transformation which takes
the standard basis into $\{\boldsymbol{\al},\,\boldsymbol{\b},\,\boldsymbol{\g}\}$.

\section{Polar Decomposition }

One can try to mimic the classical polar representation of complex
numbers in the current 3D situation as follows.

Let $t=|\boldsymbol{\mathbf{\mathit{\boldsymbol{T}}}}|$ be the modulus
of $\boldsymbol{\mathbf{\mathit{\boldsymbol{T}}}}=u+\j v+\j\j w\!\in\!\R^{3}$,
and let $\alpha,\beta,\gamma$ be angles between the corresponding
vector $\overrightarrow{T}$ and the positive direction of coordinate
axes $e_{1},e_{2}$ and $e_{3}$, respectively. Then
\[
cos\,\alpha=\frac{u}{t},\, cos\,\beta=\frac{v}{t},\, cos\,\gamma=\frac{w}{t},\ \tag{3.1}
\]
and we may write down $\boldsymbol{\mathbf{\mathit{\boldsymbol{T}}}}\,$
in the form: 
\[
\boldsymbol{\mathbf{\mathit{\boldsymbol{T}}}}=t\cdot(cos\,\alpha+\j cos\,\beta+\j\j cos\,\gamma),\ \tag{3.2}
\]
where $\cdot\,$ stands for the multiplication of a real scalar by
a $J_{3}$-number and angles $\alpha$, $\beta$ and $\gamma$ are
related by 

\[
cos\text{\texttwosuperior\ensuremath{\alpha}}+cos{{}^2}\beta+cos\text{\texttwosuperior\ensuremath{\gamma}}=1.\ \tag{3.3}
\]

By putting $t=1$ in (3.2) one gets the operator: 

\[
\T_{1}=cos(\alpha)+\j cos(\beta)+\j\j cos(\gamma),\ \tag{3.4}
\]
which acts on a $J_{3}$-number $\S=x+y+\j\j z$ by the following
rule:

\[
\begin{matrix}\S\d\T_{1}= & x\cdot cos(\text{\ensuremath{\alpha}})-z\cdot cos(\text{\ensuremath{\beta}})-y\cdot cos(\text{\ensuremath{\gamma}})+\\
 & \j(y\cdot cos(\text{\ensuremath{\alpha}})+x\cdot cos(\text{\ensuremath{\beta}})-z\cdot cos(\text{\ensuremath{\gamma}}))+\\
 & \j\j(z\cdot cos(\alpha)+y\cdot cos(\beta)+x\cdot cos(\gamma))
\end{matrix}\ \tag{3.5}
\]

Notice the resemblance with the rotation operation in the complex
plane:
\[
z\cdot t_{1}=(x\cdot cos(\text{\ensuremath{\theta}})-y\cdot sin(\text{\ensuremath{\theta}}))+\i\,(y\cdot cos(\text{\ensuremath{\theta}})+x\cdot sin(\text{\ensuremath{\theta}})),\ \tag{3.6}
\]

where $z=x+\i\, y$ is an arbitrary complex number and $t_{1}=cos(\theta)+\i\, sin(\theta)$
is an operator which rotates $z$ by $\theta$ radians counterclockwise
about origin.

\textbf{Definition 3.1.} We will call the $J_{3}$-number $t_{1}$
the \textbf{\textit{direction}} of $\boldsymbol{\mathbf{\mathit{\boldsymbol{T}}}}$
and denote it by $dir^{J_{3}}(\boldsymbol{\mathbf{\mathit{\boldsymbol{T}}}})$:

\[
dir^{J_{3}}(\boldsymbol{\mathbf{\mathit{\boldsymbol{T}}}})=cos\,\alpha+\j cos\,\beta+\j\j cos\,\gamma\ \tag{3.7}
\]

By recollecting (1.4) and (3.1) one can write down the representation
of $\mathbf{\mathit{\boldsymbol{\boldsymbol{T}}}}$ in the following
\textit{polar} form:

\[
\boldsymbol{\mathbf{\mathit{\boldsymbol{T}}}}=|\boldsymbol{\mathbf{\mathit{\boldsymbol{T}}}}|\cdot dir^{J_{3}}(\boldsymbol{\mathbf{\mathit{\boldsymbol{T}}}}).\ \tag{3.8}
\]

Note that though (3.2) or (3.8) formally resemble the polar representation
(a trigonometric form) of a complex number, some basic properties
of the latter fail to survive in 3D. For example, the modulus of the
$\d$-product of $J_{3}$-numbers is not always equal to the product
of moduli, see Lemma \textbf{5}. 

\smallskip{}

We are going to present another polar form of any $J_{3}$-number
$\T$:

\[
\T=\P\d\boldsymbol{U},\ \tag{3.9}
\]

where $\P,\:\boldsymbol{U\,\in\R^{3}}$ are such that $\P^{2}=\T\T^{\boldsymbol{*}}$
($\P$ plays a role of the modulus) and $\boldsymbol{UU^{*}}=\boldsymbol{1}$
($\boldsymbol{U}$ plays a role of the direction). 

\medskip{}

To this end, let us take a $J_{3}$-number $\S=a\cdot\al+b\cdot\b+c\cdot\g$
and denote the left hand side of the formula (2.14) by $\Q$:

\[
\Q=\S\d\S^{*}=a^{2}\al+(b^{2}+c^{2})\b.\ \tag{3.10}
\]

Since by \textbf{Remark} \textbf{2 }of Subsection\textbf{ 2.3 }the
altitude $|||\Q|||=a^{2}$ is non-negative, the monic equation

\[
\X^{2}=\Q\ \tag{3.11}
\]

is solvable by Theorem \textbf{1.7}. Among its possible solutions
\[
\X=\pm|a|\al+(\pm\sqrt{b^{2}+c^{2}})\b\ \tag{3.12}
\]
we will choose as a $\P$ the one with the non-negative coefficients:

\[
\boldsymbol{P}=|a|\al+\sqrt{b^{2}+c^{2}}\b.\ \tag{3.13}
\]

If $\S$ is $\d$-invertible then by Theorem 1.2 one has
\[
||\S||=a(b^{2}+c^{2})\neq0,
\]

i.e., both $a\neq0$ and $\, b^{2}+c^{2}\neq0$ in which case $\P$
is also \textit{$\d$}-invertible with its \textit{$\d$}-inverse:
\[
\boldsymbol{P}^{-1}=\frac{1}{|a|}\al+\frac{1}{\sqrt{b^{2}+c^{2}}}\b.\ \tag{3.14}
\]

Indeed, by invoking the multiplication Table \textbf{2} and the equality
(1.37) one can easily check that $\boldsymbol{P}^{-1}\d\P=\al+\b=\boldsymbol{1}$.
This enables us to divide $\S$ by $\P$ to obtain the following $J_{3}$-number:
\[
\boldsymbol{U=\S\d\P^{-1}=}a'\al+b'\b+c'\g,\ \tag{3.15}
\]

where $a'=\frac{a}{|a|},\, b'=\frac{b}{\sqrt{b^{2}+c^{2}}},\, c'=\frac{c}{\sqrt{b^{2}+c^{2}}}$.
Since $a'^{2}=1$ and $b'^{2}+c'^{2}=1$ it immediately follows from
(2.14) and (1.37) that

\[
\boldsymbol{U^{*}\d}\,\boldsymbol{U}=\al+\b=\boldsymbol{1}.\ \tag{3.16}
\]

Thus, any $\d$-invertible $J_{3}$-number $\S=a\cdot\al+b\cdot\b+c\cdot\g$
could be expressed as the following $J_{3}$-product: 

\[
\S=(|a|\al+r\b)\d(\frac{a}{|a|}\al+cos(\theta)\b+sin(\theta)\g),\ \tag{3.17}
\]

where $r=\sqrt{b^{2}+c^{2}}$ and $\theta=arcsin(\nicefrac{c}{r}$)
if $c\geq0$, or $\theta=\pi-arcsin(\nicefrac{c}{r}$) if $c<0$ (note
that both $a\neq0$ as well as $r>0$ do not vanish due to the $\d$-invertibility
of $\S$).

Let us observe that the first factor $\P=|a|\al+r\b$ which would
be referred to as the $J_{3}$-modulus belongs to the positive quadrant
of the conjugate plane $R$ while the second one $\boldsymbol{U}=\frac{a}{|a|}\al+cos(\theta)\b+sin(\theta)\g$
belongs to one of two circles parallel to $M$ and thus $\boldsymbol{U}$
defines some rotation around the $L$ axis. 

One can rewrite (3.17) in the following form for an arbitrary $\S\in\R^{3},$
including zero divisors:

\[
\S=(a\al+r\b)\d(\al+cos(\theta)\b+sin(\theta)\g),\, a\in\R,\, r\geq0,\,-\pi<\theta\leq\pi.\ \tag{3.18}
\]

This representation has an important property. Namely, if
\[
\S_{i}=(a_{i}\al+r_{i}\b)\d(\al+cos(\theta_{i})\b+sin(\theta_{i})\g),\, i=1,2,
\]

then due to commutativity of the $\d$-product, and by invoking the
multiplication \textbf{Table 2} and elementary trigonometric identities
one has

\[
\S_{1}\d\S_{2}=(a_{1}a_{2}\al+r_{1}r_{2}\b)\d(\al+cos(\theta_{1}+\theta_{2})\b+sin(\theta_{1}+\theta_{2})\g).\ \tag{3.19}
\]

Moreover, the formula (3.18) represents the entire 3D space as the
$\d$-product of a half-plane and a circle:

\[
\R^{3}=R_{+}\d C_{\boldsymbol{\alpha}},\ \tag{3.20}
\]
 where $R_{+}$ is a half-plane of the conjugate plane $R$ and $C_{\boldsymbol{\alpha}}$
is an $\al$-centered circle which is parallel to the $\j$-invariant
plane $M$. 

\smallskip{}

\textbf{Remark}. Note that in view of (3.18) one can parametrize any
$J_{3}$-number 

\[
\boldsymbol{\mathbf{\mathit{\boldsymbol{T}}}}=u+\j v+\j\j w=a\cdot\al+b\cdot\b+c\cdot\g
\]
by a triple $(r,\,\theta,\, a)$ which could be interpreted as the
\textit{cylindrical coordinates} with respect to the reference $L$
axis and $M$ plane. Here

\[
r=\sqrt{b^{2}+c^{2}}=\sqrt{|\mathit{\boldsymbol{\boldsymbol{T}}}|^{2}+\left\{ \T\right\} },\: a=u-v+w=|||\T|||\ \tag{3.21}
\]

are linear parameters (\textit{radius} and \textit{altitude,} respectively),
while \textit{azimuth} $\theta$ is an angle between the \textit{polar}
axis $\b$ and an orthogonal projection $O\T'$ of the vector $\overrightarrow{\T}$
on the reference plane $M$ (see Fig. 3 below). 

\begin{figure}
\includegraphics{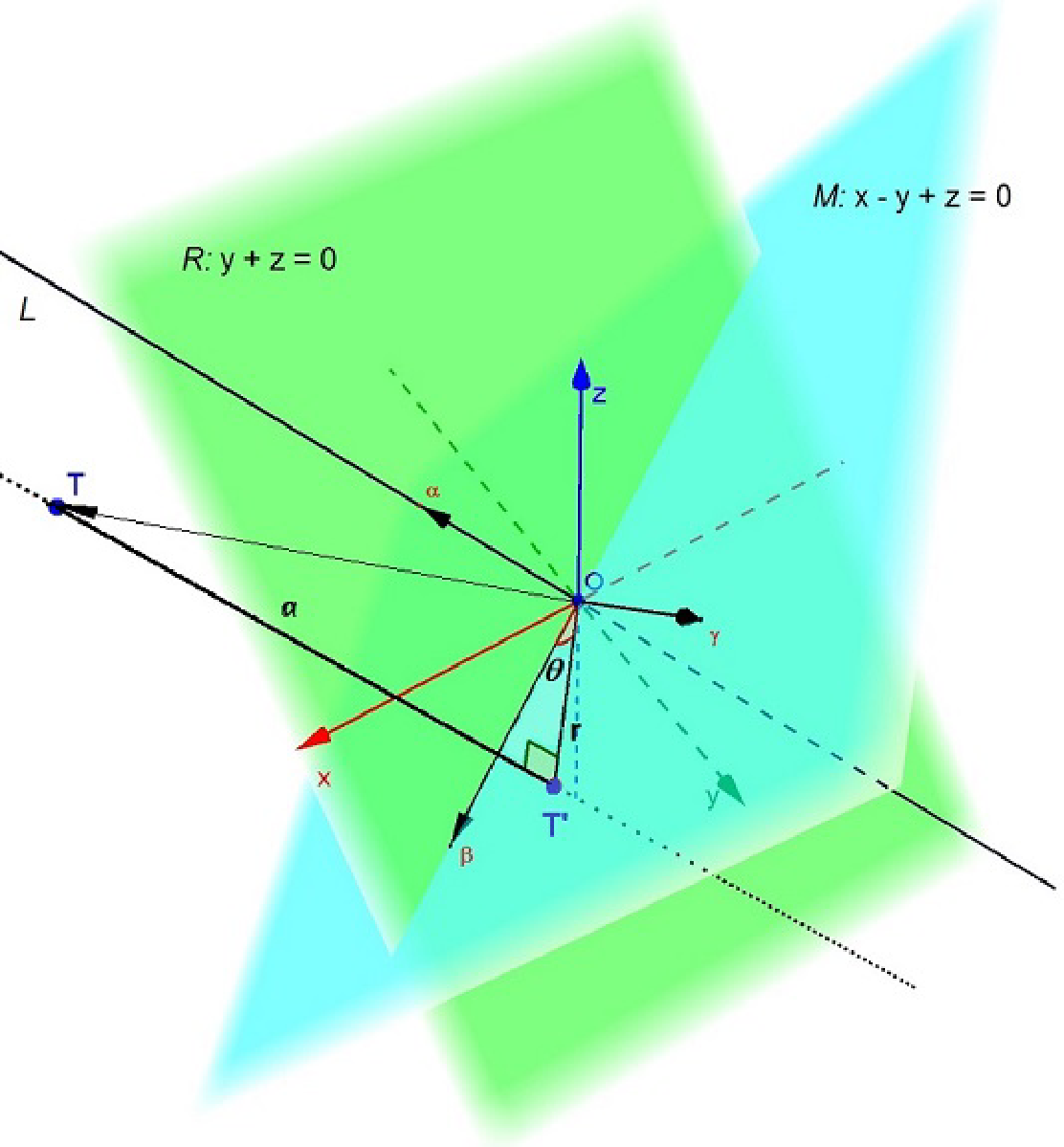}

\caption{Cylindrical coordinates (r, $\theta,\, a)$ of $\T=a\cdot\al+b\cdot\b+c\cdot\g$
with respect to $L$ and $M$ }

\end{figure}

\newpage{}

\section{Some Elementary Functions in $\R^{3}$}

In what follows we are going to demonstrate the use of $J_{3}$-numbers
as arguments of higher-order polynomials and some elementary functions,
in particular the exponential function.

\subsection{Polynomials of $J_{3}$-numbers and power series}

Since the $\d$-multiplication of $J_{3}$-numbers is a closed operation
in $\R^{3}$, one can recursively define the integer powers of a $J_{3}$-number
as follows:

\[
\X^{n}=\X^{n-1}\d\X,\:\X^{0}=\boldsymbol{1},\:\X^{1}=\X,\; n=1,2,...\ \tag{4.1}
\]

and consider polynomials in higher powers with coefficients $\boldsymbol{A}_{k}$
being $J_{3}$-numbers like it was done before for the quadratic in
Subsection \textbf{1.10}:
\[
P(\X)=\sum_{k=0}^{n}\boldsymbol{A}_{k}\d\X^{k}.\ \tag{4.2}
\]

Obviously, all the familiar algebraic manipulations with the real
or complex polynomials remain unchanged for polynomials in $\R^{3}$.
Similarly to the classical case we could also introduce formal infinite
power series:

\[
s(\X)=\sum_{k=0}^{\infty}\boldsymbol{A}_{k}\d\X^{k},\;\X\in\R^{3}.\ \tag{4.3}
\]

\subsection{$J_{3}$-trigonometric functions}

Let us remind that the classical trigonometric functions $sin$ and
$cos$ have the following power series representation

\[
cos(x)=\sum_{n=0}^{\infty}\frac{(-1)^{n}x^{2n}}{(2n)!},\; sin(x)=\sum_{n=0}^{\infty}\frac{(-1)^{n}x^{2n+1}}{(2n+1)!}.\ \tag{4.4}
\]

Motivated by Euler's formula that gives a connection between complex
numbers, exponents and trigonometry

\[
e^{\i x}=cos\, x+\i\, sin\, x,\; x\in\R\ \tag{4.5}
\]

we define three $J_{3}$-trigonometric functions of a real variable
$x$:

\[
cos_{0}(x)=\sum_{n=0}^{\infty}\frac{(-1)^{n}x^{3n}}{(3n)!},\ sin_{1}(x)=\sum_{n=0}^{\infty}\frac{(-1)^{n}x^{3n+1}}{(3n+1)!},\ sin_{2}(x)=\sum_{n=0}^{\infty}\frac{(-1)^{n}x^{3n+2}}{(3n+2)!}.\ \tag{4.6}
\]

It is easy to see that these power series are uniformly convergent
and that $cos_{0},\, sin_{1}$, and $sin_{2}$ are all smooth functions
that satisfy the following differential equations:
\[
\frac{d}{dx}sin_{2}(x)=sin_{1}(x),\,\frac{d}{dx}sin_{1}(x)=cos_{0}(x),\,\frac{d}{dx}cos_{0}(x)=-sin_{2}(x).\ \tag{4.7}
\]

Moreover, similarly to the classical $sin$ and $cos$ functions that
are solutions of a simple harmonic oscillator equation $y''(x)+y(x)=0$
each of the three $J_{3}$-trigonometric functions satisfy the following
third order linear differential equation:

\[
y'''(x)+y(x)=0.\ \tag{4.8}
\]

By solving (4.8) and taking into account (4.7) the $J_{3}$-trigonometric
functions (4.6) can be expressed in terms of the standard elementary
functions as follows:

\[
\begin{array}{c}
cos_{0}(x)=\frac{1}{3}e^{-x}(1+e^{\frac{3}{2}x}cos(\frac{\sqrt{3}}{2}x)),\\
\\
sin_{1}(x)=\frac{1}{3}e^{-x}(-1+e^{\frac{3}{2}x}(cos(\frac{\sqrt{3}}{2}x)+\sqrt{3}sin(\frac{\sqrt{3}}{2}x))),\\
\\
sin_{2}(x)=\frac{1}{3}e^{-x}(1-e^{\frac{3}{2}x}(cos(\frac{\sqrt{3}}{2}x)+\sqrt{3}sin(\frac{\sqrt{3}}{2}x))).
\end{array}\ \tag{4.9}
\]

\subsection{Exponential form of a $J_{3}$-number and Euler's identity in $\R^{3}$}

\textbf{Definition 4.1.} The \textbf{exponent} of a $J_{3}$-number
$\X$ is defined as follows:

\[
Exp(\X)=\sum_{k=0}^{\infty}\frac{\X^{k}}{k!},\;\X\in\R^{3}.\ \tag{4.10}
\]

Due to the isomorphism (2.4) between $J_{3}$-numbers and Toeplitz
matrices of the form (2.3) the above series converges absolutely and,
moreover, in view of the commutativity:

\[
Exp(\X+\Y)=Exp(\X)\d Exp(\Y),\;\forall\X,\,\forall\Y\in\R^{3}.\ \tag{4.11}
\]
\textbf{Remark}. Note that $Exp(\boldsymbol{0})=\boldsymbol{1}$ and
that $Exp(\X)$ is always $\d$-invertible with the inverse

\[
Exp(\X){}^{-1}=Exp(-\X),\;\forall\X\in\R^{3}.\ \tag{4.12}
\]

By inserting $\X=u+\j v+\j\j w$ into (4.10) and using (4.11) one
could get the following factorization:

\[
Exp(u+\j v+\j\j w)=Exp(u)\d Exp(\j v)\d Exp(\j\j w),\ \tag{4.13}
\]

where the first factor $Exp(u)=e^{u}$ is real while the remaining
factors could be expressed due to $\j^{3}=-\boldsymbol{1}$ as follows:

\[
Exp(\j\, v)=\sum_{k=0}^{\infty}\frac{(\j v)^{k}}{k!}=cos_{0}(v)+\j\, sin_{1}(v)+\j\j\, sin_{2}(v),\ \tag{4.14}
\]
\[
Exp(\j\j\, w)=\sum_{k=0}^{\infty}\frac{(\j\j w)^{k}}{k!}=cos_{0}(-w)-\j\, sin_{2}(-w)-\j\j\, sin_{1}(-w).\ \tag{4.15}
\]
By combining (4.13)-(4.15) one gets

\[
Exp(u+\j v+\j\j w)=f(u,v,w)+\j\, g(u,v,w)+\j\j\, h(u,v,w),\ \tag{4.16}
\]

where $f,\, g,\, h$ are as follows:

\[
f=e^{u}(cos_{0}(v)\cdot cos_{0}(-w)+sin_{1}(v)\cdot sin_{1}(-w)+sin_{2}(v)\cdot sin_{2}(-w)),\ \tag{4.16f}
\]
\[
g=e^{u}(sin_{1}(v)\cdot cos_{0}(-w)+sin_{2}(v)\cdot sin_{1}(-w)-cos_{0}(v)\cdot sin_{2}(-w)),\ \tag{4.16g}
\]
\[
h=e^{u}(sin_{2}(v)\cdot cos_{0}(-w)+sin_{1}(v)\cdot sin_{2}(-w)-cos_{0}(v)\cdot sin_{1}(-w)).\ \tag{4.16h}
\]

If one writes down $f,\, g,\, h$ in terms of standard elementary
functions by (4.9) then the following outstanding formula would emerge
after simplifications:

\[
\begin{array}{cc}
Exp(u+\j v+\j\j w)= & \frac{1}{3}e^{u-v+w}((1+2e^{3\varphi}\cos\theta)\ \ \ \ \ \ \,\ \ +\\
 & \j(-1+e^{3\varphi}(\cos\theta+\sqrt{3}\sin\theta))\ \ \ +\\
 & \j\j(1-e^{3\varphi}(\cos\theta-\sqrt{3}\cos\theta))\ \ \,
\end{array}\ \tag{4.17}
\]
where $\varphi=\frac{1}{2}(v-w),\,\theta=\frac{\sqrt{3}}{2}(v+w)$. 

\smallskip{}

In order to verify the above formula let us express:

\[
\X=u+\j v+\j\j wi=s\cdot\al+t\cdot\b+\theta\cdot\g,\ \tag{4.18}
\]
where the orthogonal basis $\{\boldsymbol{\al},\,\boldsymbol{\b},\,\boldsymbol{\g}\}$
is as in Subsection \textbf{1.8} and the coefficients $s,\, t,\,\theta$
are given by (1.41): 
\[
s=u-v+w,\, t=\frac{1}{2}(2u+v-w),\,\theta=\frac{\sqrt{3}}{2}(v+w).\ \tag{4.19}
\]

\smallskip{}

\textbf{Lemma 11. }For any real numbers $s$ and $t$ the following
identity holds true:

\[
Exp(s\,\al)\d Exp(t\,\b)=e^{s}\al+e^{t}\b.\ \tag{4.20}
\]

\textbf{Proof. }

By using $\al+\b=\boldsymbol{1}$ and the idempotent property of $\al$
and $\b$ let us evaluate each factor separately:

\[
Exp(s\,\al)=\boldsymbol{1}+\sum_{n=1}^{\infty}\frac{(s\al)^{n}}{n!}=\al+\b+\al\cdot\sum_{n=1}^{\infty}\frac{s^{n}}{n!}=\b+\al+\al(e^{s}-1)=e^{s}\al+\b,\ \tag{4.21}
\]

\[
Exp(t\,\b)=\boldsymbol{1}+\sum_{n=1}^{\infty}\frac{(t\b)^{n}}{n!}=\al+\b+\b\cdot\sum_{n=1}^{\infty}\frac{t^{n}}{n!}=\al+\b+\b(e^{t}-1)=\al+e^{t}\b,\ \tag{4.22}
\]

Finally, due to the $\d$-orthogonality of $\al$ and $\b$ one can
see:

\[
Exp(s\,\al)\d Exp(t\,\b)=(e^{s}\al+\b)\d(\al+e^{t}\b)=e^{s}\al+e^{t}\b,\ \tag{4.23}
\]
which completes the proof.\textit{\scriptsize{}$\text{\textifsymbol[ifgeo]{32}}$}{\scriptsize \par}

Note that in view of (4.11) the above lemma\textbf{ }immediately implies:

\[
Exp(s\,\al+t\,\b)=e^{s}\al+e^{t}\b.\ \tag{4.24}
\]

\smallskip{}

\textbf{Lemma 12. }For any real number $\theta$ $ $ the following
identity holds true:

\[
Exp(\theta\,\g)=\al+cos\,(\theta)\,\b+sin\,(\theta)\,\g.\ \tag{4.25}
\]

\textbf{Proof. }

By invoking (1.40) it is easy to check that $\g^{2k}=(-1)^{k}\!\b$
and $\g^{2k+1}=(-1)^{k}\b\d\!\g=(-1)^{k}\g,\, k=1,2,..$ Thus

\[
Exp(\theta\g)=\sum_{n=0}^{\infty}\frac{(\theta\g)^{n}}{n!}=\boldsymbol{1}+\sum_{k=2}^{\infty}\frac{(\theta\g)^{2k}}{(2k)!}+\sum_{k=0}^{\infty}\frac{(\theta\g)^{2k+1}}{(2k+1)!}
\]

\[
=\al+(\b+\b\sum_{k=2}^{\infty}\frac{(-1)^{k}\theta{}^{2k}}{(2k)!})+\g\sum_{k=0}^{\infty}\frac{(-1)^{k}\theta{}^{2k+1}}{(2k+1)!}\ \tag{4.26}
\]
\[
=\al+cos\,(\theta)\,\b+sin\,(\theta)\,\g.\text{\textifsymbol[ifgeo]{32}}
\]

\textbf{Remark 1.} One can regard (4.25) as the 3D analogue of the
famous Euler's formula $e^{i\varphi}=cos\,\varphi+\i\, sin\,\varphi$.

\smallskip{}

\textbf{Theorem 4.1. }For any real numbers $s,t,\, p$ the following
identity holds true:

\[
Exp(s\al+t\b+\theta\g)=e^{s}(\al+e^{t-s}cos\,(\theta)\,\b+e^{t-s}sin\,(\theta)\,\g).\ \tag{4.27}
\]

\textbf{Proof. }It is enough to invoke (4.11) and apply \textbf{Lemmas}
\textbf{11} and \textbf{12}:

\[
Exp(s\al+t\b+\theta\g)=Exp(s\al+t\b)\d Exp(\theta\g)=e^{s}\al+e^{t}cos\,(\theta)\,\b+e^{t}sin\,(\theta)\,\g,
\]

where the multiplication Table \textbf{2} has to be used.\textit{\scriptsize{}$\text{\textifsymbol[ifgeo]{32}}$}{\scriptsize \par}

\smallskip{}

\textbf{Remark 2.} Let us note that (4.27) resembles the important
classical formula $e^{x+iy}=e^{x}(cos\, y+\i\, sin\, y).$

\smallskip{}

\textbf{Theorem 4.2. (Exponential form of a $J_{3}$-number). }

If $\S=a\cdot\al+b\cdot\b+c\cdot\g$ is $\d$-invertible, then it
admits the following representation:

\[
\S=(a\al+r\b)\d Exp(\theta\g),\:\ \tag{4.28}
\]

where $r=\sqrt{b^{2}+c^{2}}$ and $\theta=arctan(\frac{b}{c}$) if
$c\neq0$, or $\theta=arccot(\frac{c}{b}$) if $b\neq0$. 

\textbf{Proof.} 

Since $\S$ is $\d$-invertible at least one of coefficients $b$
or $c$ is nonzero, and hence one can define $\theta$. Furthermore,
due to (3.18) one has: 

\[
\S=(a\al+r\b)\d(\al+cos(\theta)\b+sin(\theta)\g).\ \tag{4.29}
\]

It is enough now to combine (4.25) and (4.29).\textit{\scriptsize{}$\text{\textifsymbol[ifgeo]{32}}$} 

\smallskip{}

When plugging $\theta=2\pi k$ into (4.25) one gets the following
formula:

\[
Exp(2\pi k\g)=\boldsymbol{1},\, k\in\mathbb{Z}\ \tag{4.30}
\]

which resembles the famous Euler's identity: $e^{\i(\pi+2\pi k)}=-1,\, k\in\mathbb{Z}\text{.}$

\smallskip{}

\textbf{Remark 3. }There is no $J_{3}$-number $\X$ such that $Exp(\X)=\boldsymbol{-1}$.
This follows from the fact that according to (4.27) the altitude $|||Exp(\X)|||=e^{s}>0$
while $|||\boldsymbol{-1}|||=-1$.

\smallskip{}

\textbf{Theorem 4.3. }The exponential curve $Exp(\theta\g),\:\theta\in\mathbb{R},$
is a circle in $\R^{3}$.

\textbf{Proof. }Let us observe that all the points of the curve $Exp(\theta\,\g)=\al+cos\,\theta\,\b+sin\,\theta\,\g$
are at the unit distance from both the origin $\boldsymbol{O}=(0,0,0)$$\,$
as well as from the $J_{3}$-invariant plane $M$. Indeed, it follows
from (1.43) and (1.60) that for $\forall\theta\in\mathbb{R}$:

\[
|Exp(\theta\g)|=\sqrt{\frac{1^{2}+2cos^{2}\theta+2sin^{2}\theta}{3}}=1,\quad|||Exp(\theta\g)|||=1.\ \tag{4.31}
\]

The unit altitude (see \textbf{Remark 2} after \textbf{Definition
7} of Subsection \textbf{1.10}) means that all the points $Exp(\theta\g)$
belong to the plane $x-y+z=1$ which is perpendicular to the axis
\textit{L }and intersects it at $\al=\frac{1}{3}(1,-1,1)$. This plane
intersects the unit sphere along a circle centered at $\al$ with
radius

\[
\boldsymbol{r}=\sqrt{1-|\al|^{2}}=\sqrt{1-\left(\frac{\sqrt{3}}{3}\right)^{2}}=\sqrt{\frac{2}{3}}=\frac{\sqrt{6}}{3},\ \tag{4.32}
\]
 (see also Fig. 3 below). The rest is plain.\textit{\scriptsize{}$\text{\textifsymbol[ifgeo]{32}}$}{\scriptsize \par}

\begin{figure}
\includegraphics{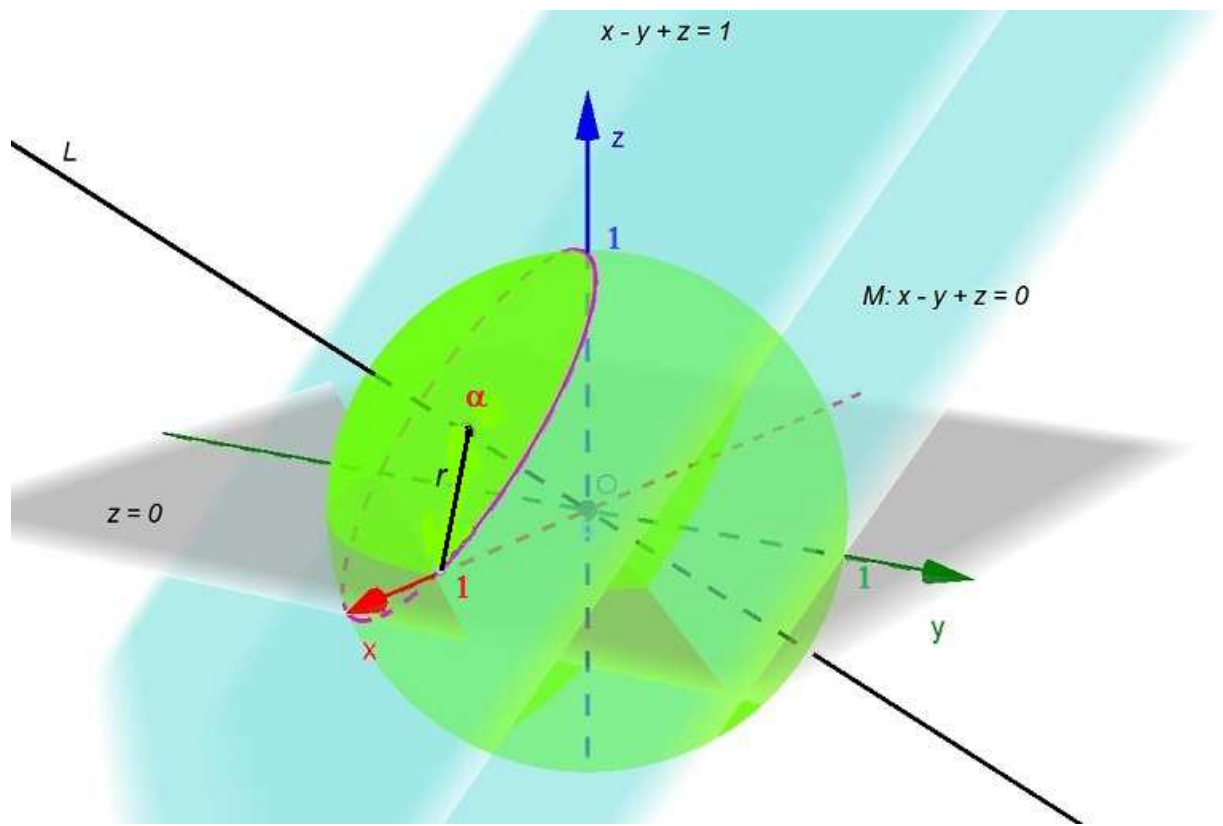}

\caption{The exponential circle (in magenta) centered at $\al$ with radius
$\boldsymbol{r}$.}
\end{figure}

\newpage{}

Finally, in order to verify formula (4.17) it is enough now to calculate
the right hand side of (4.27) according to (4.19) and (1.42) :

\[
Exp(u+\j v+\j\j w)=\frac{a+2b}{3}+\j\frac{-a+b+c\cdot\sqrt{3}}{3}+\j\j\frac{a-b+c\cdot\sqrt{3}}{3},\ \tag{4.33}
\]

where $a=e^{u-v+w},\, b=e^{\frac{1}{2}(2u+v-w)}cos(\frac{\sqrt{3}}{2}(v+w)),\, c=e^{\frac{1}{2}(2u+v-w)}sin(\frac{\sqrt{3}}{2}(v+w)).$

\subsection{Logarithm of a $J_{3}$-number}

The exponentiation of a $J_{3}$-number yields a $J_{3}$-number and
thus is a closed operation. 

It follows from (4.27) that for any $J_{3}$-number $\X=s\al+t\b+\theta\g\,$
both the altitude of $Exp(\X)$:

\[
|||Exp(\X)|||=e^{|||\X|||}=e^{s}>0,\ \tag{4.34}
\]

as well as the modulus of its projection on the plane $M$ are positive:
\[
|e^{t}cos\,(\theta)\,\b+e^{t}sin\,(\theta)\,\g|=e^{t}>0.\ \tag{4.35}
\]

This means that the exponential function in $\R^{3}$ sends any $J_{3}$-number
to $\R_{+}^{3}=\left\{ \X=x+\j y+\j\j z:\: x-y+z>0,\:\X\notin L\right\} $
- the half-space which is above the $\j$-invariant plane $M$, with
an exclusion of the $\j$-invariant line $L$. One could introduce
the logarithm function on this domain $\R_{+}^{3}$ much in the same
way as in the complex plane, namely, as an inverse to the exponential
function $Exp.$ 

\textbf{Definition 4.2.} Given a $J_{3}$-number $\Y\in\R_{+}^{3}$,
a $J_{3}$-number $\X$ is called the \textbf{logarithm} of $\Y$
(denoted by $\X=Log(\Y)$) if $\Y=exp(\X)$.

\smallskip{}

The $Log$ function in $\R^{3}$ is multivalued similarly to the classical
logarithm in the complex plane. For example, due to (4.30) one has
$Log(\boldsymbol{1})=\boldsymbol{0}+2\pi k\g,\, k\in\mathbb{Z}.$

It follows from (4.11) that the $Log$ function is subject to the
following logarithmic identity, up to $2\pi k\g$ for some $k\in\mathbb{Z}$:

\[
Log(\X\d\Y)=Log(\X)+Log(\Y),\:\forall\X,\,\Y\in\R_{+}^{3}\ \tag{4.36}
\]

If $\X\in\R_{+}^{3}$ and $n\in\mathbb{N}$ then (4.35) implies:
\[
Log(\X^{n})=n\cdot Log(\X)\ \tag{4.37}
\]

\smallskip{}

Since $|||\j|||=0-1+0=-1,\: Log(j)$ is not defined. Let us compute
now $Log(\j^{2})$ instead. It is easy to see that $|\j^{2}|=1$ and
$|||\j^{2}|||=1$, i.e., a point in $\R^{3}$, which corresponds to
the $J_{3}$-number $\j\j$, belongs to the circle (4.29), and thus
$Log(\j^{2})=\theta\g$ for some $\theta\in\mathbb{R}$. Since by
(1.41)
\[
\j^{2}=0+0\j+\j\j=\al-\frac{1}{2}\b+\frac{\sqrt{3}}{2}\g\ \tag{4.38}
\]

it follows from (4.29) that $cos\,\theta=-\frac{1}{2},\: sin\,\theta=\frac{\sqrt{3}}{2}$.
By solving this elementary trigonometric system one gets

\[
Log(\j^{2})=\frac{2\pi}{3}\g+2\pi k\g,\, k\in\mathbb{Z\ \tag{4.39}}
\]

Finally, since by the basic property (1.2) $-\j=\j^{2}\d\j^{2}$ it
is easy to compute: $Log(-\j)=Log(\j^{2}\d\j^{2})=Log(\j^{2})+Log(\j^{2})=(\frac{4\pi}{3}+2\pi k)\g$,
or

\[
Log(-\j)=(\frac{\pi}{3}+\pi(2k+1))\g,\, k\in\mathbb{Z\ \tag{4.39}}
\]

\section*{Epilogue }

The three-dimensional hypercomplex $J_{3}$-number has been introduced.
It is a scalar composed of three components which represents a point
in $\R^{3},$ similarly to the complex number representing a point
in the complex plane. The algebraic and geometric properties of $J_{3}$-numbers
have been presented. 

The analytic properties are out of the scope of the current work.
We just note that $J_{3}$-analytic $\R^{3}$-valued function $F(z)=f(u,v,w)+\j\, g(u,v,w)+\j\j\, h(u,v,w)$
of a $J_{3}$-argument $z=u+\j v+\j\j w$ could be defined for which
the analogue of the classical Cauchy - Riemann equations holds true:

\[
\frac{\partial f}{\partial u}=\frac{\partial g}{\partial v}=\frac{\partial h}{\partial w},\ \frac{\partial g}{\partial u}=\frac{\partial h}{\partial v}=-\frac{\partial f}{\partial w},\ \frac{\partial h}{\partial u}=-\frac{\partial f}{\partial v}=-\frac{\partial g}{\partial w}.
\]

Being a scalar, the $J_{3}$-numbers possess all attributes we demand
from a scalar to have, that is being associative, commutative and
distributive under addition as well as under multiplication. The reality
that $J_{3}$-numbers are scalars enable their use as arguments in
elementary functions as has been demonstrated. 

The beauty of having a number, not a vector, representing a point
in a 3D space is enchanting.

\section*{Acknowledgments}

I wish to express my gratitude to Dr. Michael Shmoish (Technion, Haifa)
for reviewing and discussing the manuscript, and for making valuable
remarks.

\section*{Bibliography}

W. R. Hamilton, On the geometrical interpretation of some results
obtained by calculation with biquaternions, Proceedings of the Royal
Irish Academy, vol. 5, pp. 388\textendash{}90, 1853.

B. L. van der Waerden, Modern Algebra, F. Ungar, New York; 3rd Edition,
1950.

B. L. van der Waerden, A History of Algebra: from al-Khwarizmi to
Emmy Noether, Springer-Verlag, Berlin, 1985.

G. E. Hay, Vector and tensor analysis, Dover Publications, 1953.

Ruel V. Churchill, Complex Variables and Applications, McGraw-Hill
Inc., US; 2nd edition, December 1960.

E. Dale Martin, A system of three-dimensional complex variables, NASA
technical report, 1986.

I. Kantor and A. Solodovnikov, Hypercomplex numbers, Springer-Verlag,
New York, 1989.

P. Kelly, R. L Panton, and E Dale Martin, Three-dimensional potential
flows from functions of a 3D complex variable, Fluid Dynamics Research
6:119-137, 1990.

G. B. Price, An Introduction to Multicomplex Spaces and Functions,
Marcel Dekker, New York, 1991.

G.Turk and M. Levoy, Zippered Polygon Meshes from Range Images, in
Computer Graphics Proceedings, ACM SIGGRAPH, pp. 311-318, 1994.

C. M. Davenport, A commutative hypercomplex algebra with associated
function theory, in Clifford Algebras With Numeric and Symbolic Computations,
pp 213-227, 1996.

S. Olariu, Complex numbers in three dimensions, arXiv:math.CV/0008120,
2000.

S. Olariu, Complex numbers in N dimensions, Elsevier, 2002.

Jian-Jun Shu, and Li Shan Ouw, Pairwise alignment of the DNA sequence
using hypercomplex number representation, Bulletin of Mathematical
Biology, Vol. 66, No. 5, pp. 1423-1438, 2004.

Wolfram Research, Inc., Mathematica, Version 5.1, Champaign, IL, 2004.

Alfsmann, D., Gockler, H.G., Sangwine, S.J., Ell, T.A., Hypercomplex
algebras in digital signal processing: Benefits and drawbacks. In:
Proc. 15th European Signal Processing Conference, pp. 1322\textendash{}1326,
2007

Anderson, M., Katz V.J., and Wilson R.J. Who Gave You the Epsilon?:
And Other Tales of Mathematical History. Washington, DC: Mathematical
Association of America, 2009.

GeoGebra. Version 5.0. URL http://http://www.geogebra.org/, 2014. 

R Core Team. R: A language and environment for statistical computing.
Version 3.2.1, R Foundation for Statistical Computing, Vienna, Austria.
URL http://www.R-project.org/, 2015.

Daniel Adler, Duncan Murdoch and others. rgl: 3D Visualization Using
OpenGL. R package version 0.95.1260/r1260. http://R-Forge.R-project.org/projects/rgl/,
2015.

Ya. O. Kalinovsky, Yu. E. Boyarinova, I. V. Khitsko, Reversible Digital
Filters Total Parametric Sensitivity Optimization using Non-canonical
Hypercomplex Number System, arXiv:cs.NA/1506.01701, 2015.

\newpage{}

\noindent \begin{center}
{\large{}Michael Shmoish}
\par\end{center}{\large \par}

{\LARGE{}Afterword }{\LARGE \par}

The story of a quest for a proper three-dimensional analogue of the
complex numbers is rich and fascinating. You probably remember the
famous question Hamilton's sons used to ask him every morning in early
October 1843: \textquotedbl{}Well, Papa, can you multiply triples?\textquotedbl{}
Sir W. R. Hamilton, according to his own letter, was always obliged
to reply, with a sad shake of the head: \textquotedbl{}No, I can only
add and subtract them.\textquotedbl{} Soon after he saw a way to multiply
quadruples leading to his prominent discovery of \textit{quaternions}.

I first met Mr. Shlomo Jacobi under sad circumstances in early 2012,
when he had just lost his beloved wife to a deadly disease. Still
Shlomo was strong enough to talk about his idea of three-dimensional
hypercomplex numbers and to show me the following ingenious multiplication
of triples 

\[
\begin{matrix}(a,\, b,\, c)\d\,(u,\, v,\, w)= & (au-bw-cv,\, av+bu-cw,\, aw+bv+cu),\end{matrix}
\]

that he discovered in early 1960s, shortly before his graduation from
the Technion.

It was extremely important to Hamilton that the modulus of a product
of two vectors would be equal to the product of their moduli. This
\textit{law of moduli} requirement (which is impossible to achieve
in dimension three due to the well-known theorems by Frobenius and
Hurwitz on real division algebras) was abandoned by Shlomo in favor
of commutativity of the above \textbf{$\boldsymbol{\d}$}-product
even though zero divisors appeared. \textquotedbl{}They only add interest\textquotedbl{}
as Olga Taussky Todd put it once. Previous attempts to introduce hypercomplex
numbers were mostly algebraic, Hamilton and his successors were trying
to devise a \textquotedbl{}wise\textquotedbl{} multiplication table.
The above article suggests a purely geometric approach by defining
a linear operator $\j$ which transforms the three-dimensional Euclidean
space into itself:

\[
\j:\,(x,\, y,\, z)\rightarrow(-z,\, x,\, y)
\]

and thus mimics the multiplicative action of imaginary unit in the
complex plane: $\,\i\cdot\,(x+\i y)=(-y+\i x)$. The $\d$-product
emerges naturally from the basic properties of operator $\j$ and
a definition of three-component \textbf{$\boldsymbol{J_{3}}$}-numbers,
while the law of moduli happens to be replaced by \textit{moduli inequality}:

\textbf{
\[
|\boldsymbol{\boldsymbol{S}}\d\mathit{\boldsymbol{\boldsymbol{T}}}|\leq\sqrt{3}\cdot|\boldsymbol{\boldsymbol{S}}|\,|\mathit{\boldsymbol{\boldsymbol{T}}}|.
\]
}

Though Shlomo's article is mainly concerned with geometric and algebraic
aspects of the $J_{3}$-numbers, many analytic properties of complex
numbers and complex-valued functions could be extended properly to
the three-dimensional case due to the above inequality.

After the discovery of quaternions some generalizations of the classical
complex numbers to higher (usually >= 4) dimensions were developed,
such as matrices, general hypercomplex number systems, and Clifford
geometric algebras. As for dimension three, I would mention an article
by Silviu Olariu where the geometric, algebraic, and analytical properties
of his tricomplex numbers, the close relatives of Shlomo's $J_{3}$-numbers,
were studied in detail. Note that Shlomo was unaware of Olariu's works
as well as earlier NASA reports by E. Dale Martin (on the theory of
three-component numbers and their applications to potential flows)
listed in the above bibliography section. I've compiled this short
bibliography which provides only a limited overview of hypercomplex-related
field and contains several mathematical textbooks from Shlomo's bookshelf. 

Shlomo's main idea was that $J_{3}$-numbers are \textit{scalars}
that could be dealt with conveniently once accustomed. Algebra $\R_{\circledast}^{3}$
of $J_{3}$-numbers is linked to geometry in three dimensions in a
simple and natural way. Based on Shlomo's mostly elementary article,
the advanced notions of invariant subspaces, idempotents, structured
matrices, algebra isomorphism could be explained easily to undergraduate
students via visualization in 3D space. The commutativity and useful
analytical properties hopefully makes the $J_{3}$-numbers a valuable
addition to the current toolbox of rotation matrices and quaternions
for manipulating objects in 3D space, optimal tracking and robotic
applications. There is also some evidence that the hypercomplex systems
similar to $J_{3}$-numbers prove useful in cryptography, physics,
digital signal processing, alignment of DNA sequences, and study of
the 3D structure of macromolecules. 

The representation of entire $\R^{3}$ as the $\d$-product of a half-plane
and a circle according to formula (3.20) might be advantageous. In
particular, one can decompose any three-dimensional body like the
Stanford bunny (in blue) into the planar part (in cyan) and the arc
(in green) and then manipulate (e.g., cluster or encrypt) each component
separately. The \textbf{$\boldsymbol{\d}$}-product would give then
a fast and easy way to recover the modified 3D object. 

\includegraphics{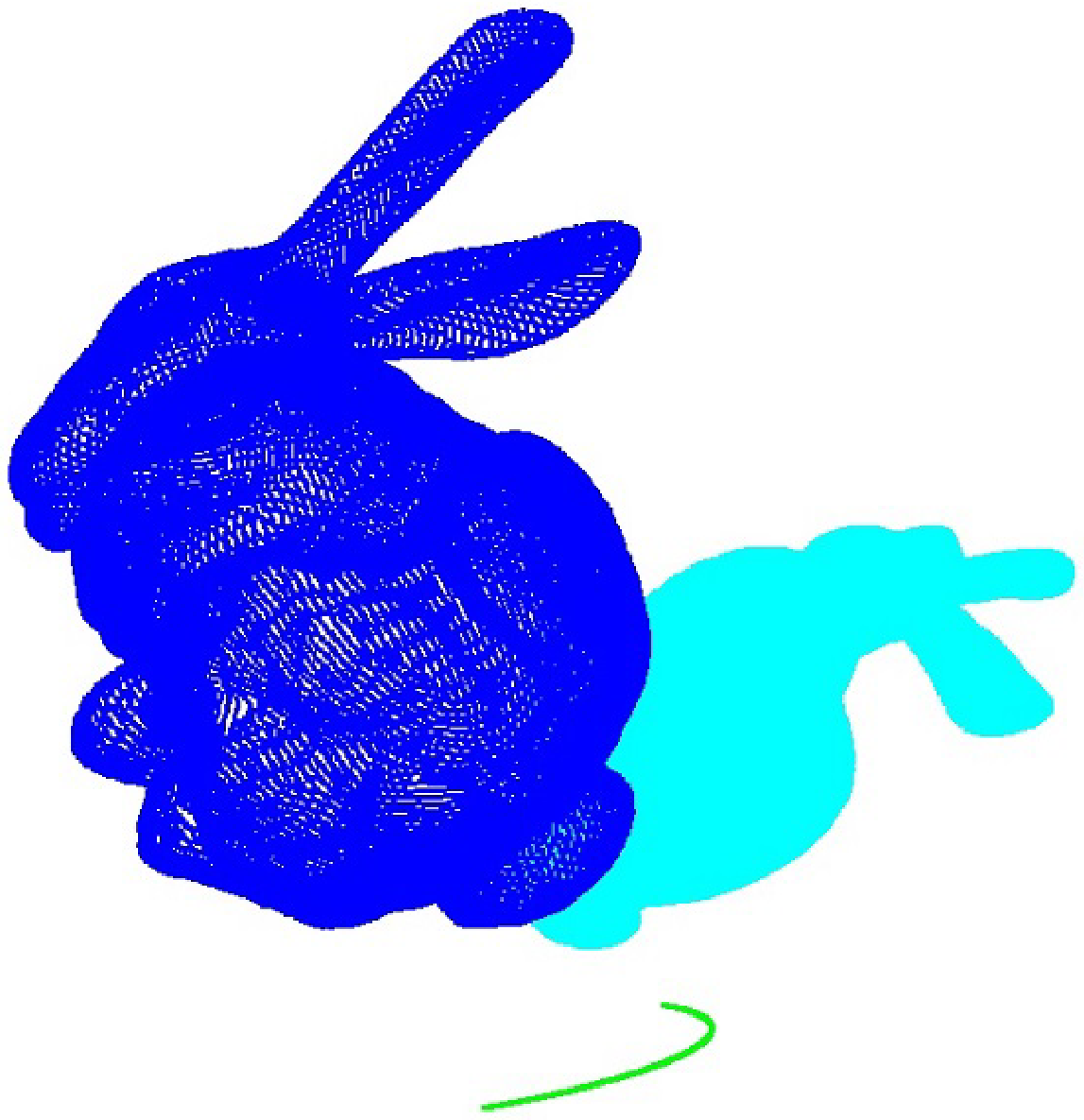}

The \textquotedbl{}bunny\textquotedbl{} illustration has been produced
using R-package 'rgl' and my R-script based on the original Shlomo's
code written in Wolfram Mathematica, while all the figures in the
above article have been produced by myself using the 3D GeoGebra.

I hope that this article will serve as a tribute to a dear friend
Shlomo Jacobi and to his life-long passion for mathematics.
\end{document}